\documentclass[12pt]{article}
\usepackage[american]{babel}
\usepackage{amsmath}
\usepackage{latexsym}
\usepackage{amssymb}
\usepackage{theorem}
\usepackage{diagrams}
\usepackage{graphicx}
\theoremstyle{change}
\newtheorem{Thm}{Theorem}[section]
\newtheorem{Cor}[Thm]{Corollary}
\newtheorem{Prop}[Thm]{Proposition}
\newtheorem{Lem}[Thm]{Lemma}
{\theorembodyfont{\rmfamily}
\newtheorem{Num}[Thm]{}

\newtheorem{Def}[Thm]{Definition}}
\newcommand{\mas}{\mathsf m}
\renewcommand{\phi}{\varphi}
\renewcommand{\rho}{\varrho}

\newcommand{\bra}[1]{\langle#1\rangle}
\newcommand{\proof}{\par\medskip\rm\emph{Proof. }}

\newcommand{\qed}{\ \hglue 0pt plus 1filll $\Box$}

\newcommand{\mapstoo}{\longmapsto}
\newcommand{\lquot}[1]{{\raisebox{-1.5pt}{$#1$}\backslash}}
\newcommand{\Herm}{\mathrm{Herm}}
\newcommand{\U}{\mathbf U}
\newcommand{\RR}{\mathbb{R}}
\newcommand{\ZZ}{\mathbb{Z}}

\newcommand{\CC}{\mathbb{C}}
\newcommand{\id}{\mathrm{id}}
\newcommand{\SKIP}[1]{}
\newcommand{\SO}{\mathbf{SO}}
\newcommand{\SL}{\mathbf{SL}}
\newcommand{\Sp}{\mathbf{Sp}}
\newcommand{\GL}{\mathbf{GL}}
\newcommand{\EU}{\mathbf{EU}}
\newcommand{\SU}{\mathbf{SU}}

\newcommand{\eps}{\varepsilon}
\renewcommand{\emptyset}{\varnothing}

\newcommand{\chr}{\mathrm{char}}

\newcommand{\Hom}{\mathrm{Hom}}

\newcommand{\Cen}{\mathrm{Cen}}
\newcommand{\hyp}{\mathrm{hyp}}

\newcommand{\cG}{\mathcal G}

\newcommand{\cO}{\mathcal O}
\newcommand{\cL}{\mathcal L}

\newcommand{\x}{\mathbf x}
\newcommand{\y}{\mathbf y}

\DeclareMathAlphabet{\varcal}{U}{rsfs}{m}{it}
\begin{document}

\title{{\bf A Maslov cocycle for unitary groups}}
\author{Linus Kramer and Katrin Tent%
\thanks{The authors were supported by SFB 478 and SFB 701}}
%\date{}
\maketitle

\begin{abstract}
We introduce a $2$-cocycle for symplectic and skew-hermitian
hyperbolic groups over arbitrary fields and skew fields, with
values in the Witt group of hermitian forms.
This cocycle has good functorial properties: it is natural under
extension of scalars and stable, so it can be viewed
as a universal $2$-dimensional characteristic class for these
groups. Over $\RR$ and $\CC$, it coincides with the first
Chern class.
\end{abstract}


\section*{Introduction}

We introduce a Maslov index and Maslov cocycle for symplectic and
hyperbolic unitary groups over arbitrary fields and skewfields.
In the classical work of Lion-Vergne \cite{LV}, this is done 
by associating to triples $(X,Y,Z)$ of
Lagrangians in a real symplectic vector space $M$ a certain integral
invariant, the Maslov index. This invariant is used to construct a
$\ZZ$-valued cocycle for the symplectic group.
The corresponding group extension of the symplectic group
is the topological
universal covering group of $\Sp_{2n}\RR$.

In this approach, it is somewhat cumbersome that one has to deal with
arbitrary triples of Lagrangians. Our starting point was the idea that
the whole construction should also work if one considers only triples
of Lagrangians in 'general position', that is, triples 
$(X,Y,Z)$ in $M$ which are pairwise opposite:
\[
M=X+Y=Y+Z=Z+X.
\]
Geometrically, such triples are much easier to classify. Moreover,
these triples carry an interesting algebraic structure.
To each pair $(X,Y)$ of opposite Lagrangians one can associate a
linear map $[Y;X]$ which identifies $X$ with the dual of $Y$ and
the dual of $X$ with $Y$.
In this way we obtain a graph, the opposition graph,
whose vertices are the Lagrangians and
whose edges join opposite Lagrangians.
Concatenating the linear maps $[Y;X]$
along closed paths in this graph, we arrive at an interesting
groupoid $\mathcal GM$, the projectivity groupoid.
A minimal closed path has length $3$,
and the resulting element in the holonomy group turns out to be
a complete geometric invariant for the triple consisting of
the three Lagrangians along the path.
This makes sense and works not just for symplectic forms, but for
arbitrary hyperbolic skew-hermitian forms over fields or
skewfields.

In order to relate this invariant to group cohomology, we need a
chain complex. A natural candidate is the flag complex of the 
opposition graph,
whose simplices are the finite complete subgraphs
(cliques). If the field is
infinite, this flag complex
is contractible, and the symplectic (or unitary) group acts on it,
so its equivariant cohomology is isomorphic to the group cohomology.

The final ingredient is the observation that along a closed path
of length $3$, the element in the holonomy group determines
a nondegenerate hermitian form, which may be viewed as an element
in a Witt group.
In this way we associate to every triangle in the opposition graph an
element in
the Witt group of hermitian forms. We verify that this map is indeed
an invariant cocycle, which gives us a $2$-cocycle for
the unitary group.

This cocycle, which we call the Maslov cocycle, 
has good functorial properties. It is stable under
direct sums of hermitian spaces and well-behaved under extension of
scalars. Furthermore, it coincides in the symplectic setting over
fields of characteristic $\neq2$ with the classical Maslov
cocycle. Our cocycle, however, exists over arbitrary fields and
skewfields of any characteristic. Furthermore, the cocycle can
be reduced to a subgroup of the Witt group, the kernel of the
signed discriminant.

The classical Maslov cocycle is important, as it yields a central
extension of the symplectic group. The question which extension is
defined by our general Maslov cocycle can by and large be reduced
to a map in algebraic $K$-theory. In the smallest case
$\Sp_2D=\SL_2D$ this is due to
Nekovar \cite{Ne} and Barge \cite{Ba}.
But even in the classical situation of a symplectic group 
$\Sp_{2n}D$ over a field $D\neq\RR$, our result appears to be the
first complete proof for this.
In general, the cocycle is related to certain symbols 
and depends on algebraic properties of the field. We carry this out
in some detail for local fields. For $\RR$ and $\CC$ the Maslov cocycle
'is' the first Chern class $c_1$ and gives the
universal covering groups of $\Sp_{2n}\RR$ and $\SU(n,n)$.
Over nonarchimedean local fields, we obtain a covering of degree
$\leq2$.

A Witt group valued Maslov cocycle appears already in \cite{LV}.
Besides this, our paper is influenced by \cite{Ne}, \cite{PPS}
(but see the remarks after \ref{CocycleTheorem}).
The idea of a 'partially defined cocycle' seems
to go back to Weil and appears also in a topological context
in \cite{Maz}. The opposition graph is used (in a different way)
in \cite{No}. The Maslov index itself has been generalized in
several ways \cite{CLM} \cite{NO}. Buildings \cite{Ti,LK}
are not mentioned in this paper, although
the motivation for our approach is the opposition relation in
spherical buildings.
Luring behind the linear algebra is the projectivity
groupoid for spherical buildings, which was first studied
systematically by Knarr \cite{Knarr} for spherical buildings of
rank $2$.

We assume that the reader is familiar with basic homological algebra,
as well as hermitian forms and unitary groups. Apart from this,
we tried to make the paper self-contained and accessible to
non-experts.

\medskip\noindent\textbf{Acknowledgement.}
Part of this work was completed while the authors were
at the School of Mathematics, Birmingham, UK.
We thank Theo Grundh\"ofer, Karl-Hermann Neeb, Chris Parker,
Andrew Ranicki, and Winfried Scharlau.

\newpage
\tableofcontents
\section{Lagrangians and hyperbolic modules}
In this section we introduce
some standard terminology from the theory of hermitian forms.
Everything we need can be found in \cite{HOM,Knus,Sch}.
We work over a field or division ring $D$ of arbitrary characteristic.
The modules we consider are finite dimensional right $D$-modules.
We assume that $J$ is an involution of $D$, i.e. an antiautomorphism
whose square is the identity (we allow $J=\id$).
The involution extends naturally to an involution
of the matrix ring $D^{n\times n}$ which we also denote by $J$.
For $\eps=\pm 1$ we put 
\[
D^\eps=\{a\in D\mid a-a^J\eps=0\}.
\]
\begin{Num}\textbf{Forms}\\
A \emph{form} on a right $D$-module $M$ is a biadditive map
$f:M\times M$ with the property that
\[
f(ua,vb)=a^Jf(u,v)b
\]
for all $u,v\in M$ and all $a,b\in D$. An 
\emph{$\eps$-hermitian form} $h$ is a form with the additional
property that 
\[
h(u,v)=h(v,u)^J\eps,
\]
and $(M,h)$ is called a \emph{hermitian module}. 
If $f$ is any form, then 
\[
h_f(u,v)=f(u,v)+f(v,u)^J\eps
\]
is $\eps$-hermitian. The hermitian forms which arise in this way are 
called \emph{trace $\eps$-hermitian} or \emph{even}. If $\chr(D)\neq 2$,
every $\eps$-hermitian form is automatically trace $\eps$-hermitian;
this is also true in characteristic $2$
if $J$ is an involution of the second kind, i.e. if
$J|_{\Cen(D)}\neq\id$, but may fail otherwise \cite[6.1.2]{HOM}.
Note also that $h_f(u,u)=0$ is equivalent to $f(u,u)\in D^{-\eps}$.
\end{Num}
\begin{Num}
The dual $M^\vee$ of $M$ (which is a left $D$-module)
can be made into a right $D$-module $M^J$
by twisting the scalar multiplication with $J$, i.e. by
setting 
\[
\xi a=[v\mapstoo a^J\xi(v)]
\]
(where $a\in D$, $\xi\in M^\vee$ and
$v\in M$). Thus forms are just linear
maps $M\rTo M^J$. A form is called \emph{non-degenerate} it the
associated linear map is injective (and hence bijective).
There is a natural notion
of an isomorphism (or isometry) of forms; the automorphism group
of a non-degenerate $\eps$-hermitian form is the 
\emph{unitary group}
\[
\U(M,h)=\U(M)=\{g\in\GL(V)\mid h(u,v)=h(g(u),g(v))\text{ for all }
u,v\in M\}.
\]
\end{Num}
\begin{Num}\textbf{Lagrangians}\\
For any subset $X\subseteq V$ we have the subspace
$X^\perp=\{u\in M\mid h(x,u)=0\text{ for all }x\in M\}$, the 
\emph{perp}. A subspace
which is contained in its own perp is called \emph{totally isotropic}
and a subspace which coincides with its perp is called a
\emph{Lagrangian}. A nondegenerate hermitian form which admits
Lagrangians is called \emph{metabolic}.
\end{Num}
\begin{Num}\textbf{The hyperbolic functor}\\
Given a right $D$-module $X$, there is a natural form $f$ on
$M=X\oplus X^J$, given by $f((x,\xi),(y,\eta))=\xi(y)$. The associated
trace $\eps$-hermitian form 
\[
h_X((x,\xi),(y,\eta))=\xi(y)+\eta(x)^J\eps
\]
(and every isometric hermitian module) is called
\emph{hyperbolic}. Obviously, $X$ is a Lagrangian, so hyperbolic
modules are metabolic. The converse is true for trace valued
hermitian forms, hence in particular in characteristic $\neq 2$
\cite[I 3.7.3]{Knus}.
The \emph{rank} of a hyperbolic module is the dimension
of $X$ (i.e. half the dimension of the hyperbolic module).
We note that the assignment
\[
\hyp:X\mapstoo (X\oplus X^J,h_X)
\]
is a functor from $D$-modules to hermitian modules, and
that $\hyp$ induces an injection\linebreak
$\GL(X)\rTo\U(X\oplus X^J)$.
\end{Num}
\begin{Num}\label{ClassesOfgroups}%
\textbf{Special cases and Lie groups}\\
Every hyperbolic form $(M,h)$ can be reduced
to one of the following three types. 

\medskip
\noindent
\textbf{Symplectic groups:}
$(J,\eps)=(\id,-1)$. Then $D$ is necessarily commutative
and $\U(M)=\Sp(M)$ is the \emph{symplectic group}.
For $M=\RR^{2n},\CC^{2n}$, these Lie groups are often denoted
$\Sp(n,\RR)$ and $\Sp(n,\CC)$.

\medskip
\noindent
\textbf{Hyperbolic orthogonal groups:}
$J=\id$ and $\eps=1\neq -1$.  Then $D$ is commutative and of characteristic
different from $2$. The group
$\U(M)=\mathbf O(M)$ is the hyperbolic \emph{orthogonal} group;
for $\RR$ and $\CC$, these Lie groups are often denoted
$\mathbf O(n,n)$ and $\mathbf O(2n,\CC)$.
We will see in \ref{CasesOrbitTypes}
below that the Maslov cocycle is uninteresting in
this situation.

\medskip
\noindent
\textbf{Standard hyperbolic unitary groups:}
If $J\neq\id$ then $\U(M)$ is the \emph{standard hyperbolic unitary group}.
Scaling the hermitian form by a suitable constant and changing the
involution,
we can assume that $\eps=-1$ (''Hilbert 90'', see \cite[p.~211]{HOM}).
The $-1$-hermitian forms are also called \emph{skew hermitian}.
Examples of involutions are the standard conjugation 
$z\mapsto\bar z$ on $\mathbb C$
and on the real quaternion division algebra
$\mathbb H$. Note that there is also the 'nonstandard' involution
$z^\alpha=-i\bar zi$ on $\mathbb H$. The skew
hyperbolic unitary groups corresponding to
$(\CC^n,z\mapsto \bar z)$, $(\mathbb H^n,z\mapsto\bar z)$ and 
$(\mathbb H^n,z\mapsto z^\alpha)$ are the Lie groups denoted
$\U(n,n)$, $\SO^*(4n)$ and $\Sp(n,n)$ in \cite[X~Tab.~V]{He}.
\end{Num}

\section{The opposition graph and triples of Lagrangians}

In this section we construct an invariant $\kappa$
which classifies triples
of pairwise opposite Lagrangians in a $-\eps$-hermitian
hyperbolic module up to isometry.
The invariant is a nondegenerate $\eps$-hermitian form.
In particular, we will have to work simultaneously with $\eps$- and
$-\eps$-hermitian forms.
We assume throughout that $M$ is a $-\eps$-hermitian hyperbolic module 
and we let 
\[
\cL=\cL(M)=\{X\in M\mid X=X^\perp\}
\]
denote its set of Lagrangians.
\begin{Def}
We call two Lagrangians $X$ and $Y$ \emph{opposite} if $X\cap Y=0$
or, equivalently, if $M=X+Y$.
If the rank of $M$ is $1$, then Lagrangians are $1$-dimensional,
and $X$ is opposite $Y$ if and only if $X\neq Y$.
\end{Def}
\begin{Lem}
\label{Card1}
If $M$ has rank $1$, then $\cL$ has $|D^\eps|+1$ elements.

\proof Let $x$ be a nonzero
vector in the $1$-dimensional space
$X$ and let $\xi\in X^J$ be its dual, i.e. $\xi(x)=1$.
Then $x$ and $\xi$ span $X\oplus X^J\cong M$.
The vector $v=(xa,\xi)$ spans a Lagrangian
if and only if $\xi(xa)=a\in D^\eps$.
There is precisely one additional 
Lagrangian, spanned by $(x,0)$.
\qed
\end{Lem}
Later it will be important that there are enough Lagrangians. We note that
$D^\eps$ is infinite if $D$ is an infinite field, unless $J=\id$ and
$\eps=-1\neq 1$. If $D$ is not commutative, then $D^\eps$ is always
infinite \cite[6.1.3]{HOM}.
\begin{Prop}
\label{EnoughOpposites}
If $|D^\eps|\geq k$, then there exists for every
finite collection $X_1,\ldots,X_k$ of Lagrangians
a Lagrangian $Y$ opposite to
$X_1,\ldots,X_k$. 

\proof 
Let $n$ denote the rank of $M$. We proceed by
induction on $k\geq1$, modifying the proof in \cite[3.30]{Ti}. 
Let $X_1,\ldots,X_k$ be $k$ Lagrangians. 
We choose a Lagrangian $Y$ such that $\ell=\dim(Y\cap X_1)$
is as small as possible, and (by the inductions hypothesis)
such that $Y$ is opposite $X_2,\ldots,X_k$. We claim
that $\ell=0$. Otherwise, we can choose a subspace
$Q\subseteq X_1$ of dimension $n-1$, such that
$X_1=Q+(Y\cap X_1)$. Now $Q^\perp$ can be
split as $Q\bot H$, with $H$ hyperbolic of rank $1$.
The $1$-dimensional Lagrangians $P$ of $H$ parameterize
the Lagrangians of $M$ containing $Q$ bijectively via
$P\mapsto Q\oplus P$. Let $P_1=Y\cap H$.
For $\nu=2,\ldots,k$, each $X_k$ determines a unique
$1$-dimensional Lagrangian $P_k\subseteq H$ with
$\dim((Q+P_k)\cap X_k)\neq 0$.
By Lemma \ref{Card1} we may choose a one-dimensional
Lagrangian $P'\subseteq H$ different from $P_1,\ldots, P_k$.
Then $Y'=P'\oplus Q$ is a Lagrangian opposite
$X_2,\ldots,X_k$ with $\dim(Y'\cap X_1)=\ell-1$, a
contradiction.
\qed
\end{Prop}
In particular, there exists always a Lagrangian $Y$ opposite
a given Lagrangian $X$.
The map $y\mapsto h(y,-)|_X$ is an isomorphism
$Y\rTo^\cong X^J$ and we have thus a unique
isomorphism of hyperbolic modules
$X\oplus X^J\rTo^\cong X\oplus Y=M$ 
extending the inclusion $X\rInto M$. If
$(X',Y')$ is another such pair, then we can choose an
linear isomorphisms $X\cong X'$ and obtain isomorphisms
\[
X\oplus Y\rTo^\cong X\oplus X^J\cong X'\oplus {X'}^J
\lTo^\cong X'\oplus Y'.
\]
Hence we have established the following result
(which also follows from Witt's Theorem \cite[6.2.12]{HOM}).
\begin{Lem}
The unitary group $\U(M)$ acts transitively on ordered
pairs of opposite Lagrangians.
\qed
\end{Lem}
\begin{Num}
\label{MatrixConvention}
We now study this $\U(M)$-action in more detail. We fix
a $D$-module $X$ of dimension $n$, with basis $\x$.
We put $Y=X^J$ and we let $\y$ denote the dual
basis. Then $M=X\oplus Y$ is hyperbolic of rank $n$, with
basis $\x,\y$, and we may work with 
$2\times 2$ block matrices. The hermitian form
$h=h_X$ on $M$ is represented by the matrix
\[
h=\scriptstyle\begin{pmatrix} 0&-\eps\\ 1&0\end{pmatrix}.
\]
We find
that the $\U(M)$-stabilizer $L$ of the ordered pair $(X,Y)$
consists of matrices of the form
\[
\ell_a=\scriptstyle\begin{pmatrix} a^{-J}&0\\0&a\end{pmatrix},
\]
with $a\in\GL_nD$ and $\ell_a\ell_{a'}=\ell_{aa'}$,
while the $\U(M)$-stabilizer
$U$ of $(X,\x)$ consists of matrices of the form
\[
u_t=\scriptstyle\begin{pmatrix} 1&t\\0&1\end{pmatrix},
\]
with $t-t^J\eps=0$, i.e. $t\in D^{n\times n}$ has to be 
$\eps$-hermitian. Note also that $u_tu_{t'}=u_{t+t'}$,
$u_t^{-1}=u_{-t}$, and that
 \[
\ell_au_t\ell_a^{-1}=u_{a^{-J}ta^{-1}}.
\]
The $\U(M)$-stabilizer $P$ of $X$ splits therefore as a
semidirect product $P=LU$, with \emph{Levi factor} $L$ and
\emph{unipotent radical} $U\unlhd P$.
\end{Num}
Next, we note that if $Z$ is another Lagrangian opposite
$X$, then we have a unique isomorphism
$X\oplus Y\rTo X\oplus Z$ fixing the basis
$\x$. This isomorphism is therefore given by an
element of
the group $U$, and we have the following result.
\begin{Lem}
The group $U$ acts regularly on the set $X^{opp}$ of
all Lagrangians opposite $X$.
\qed
\end{Lem}
Let $u_t\in U$. Then
the Lagrangian $Z=u_t(Y)$ is opposite $Y$ if and only if $M$ is
spanned by $\y,u_t(\y)$.
With the matrix notations we established before, we
have 
\[\textstyle
u(\y_\nu)=\y_\nu+\sum_\mu \x_\mu t_{\mu,\nu}.
\]
A necessary and sufficient condition for $Z=u_t(Y)$ being opposite
$Y$ is thus that the matrix $t$ is invertible.
\begin{Num}
We let
$H=\{t\in D^{n\times n}\mid t-t^J\eps=0\}$ denote the set of
all $\eps$-hermitian $n\times n$-matrices. There is a natural
left action $(a,t)\mapsto a^{-J}ta^{-1}$ of $\GL_nD$ on $H$, and 
we denote the orbit of $t$ by $\bra t$. The
orbit space
\[
\Herm_\eps(n)=\left\{\bra t\mid t\in H\right\}=\lquot LH
\]
consists thus of the isomorphism
classes of $\eps$-hermitian forms on $D^n$. We denote the subset
corresponding to the nonsingular hermitian forms by 
$\Herm_\eps^\circ(n)$.
Then we have an $L$-equivariant bijection
\[
H\rTo X^{opp}\qquad
t\mapstoo u_t(Y).
\]
Factoring out the $L$-action, we get bijections
\[
\Herm_\eps(n)\rTo \lquot L X^{opp}\quad\text{ and }\quad
\Herm_\eps^\circ(n)\rTo \lquot L(X^{opp}\cap Y^{opp})
\]
While the isomorphism $H\rTo U$ depends on the chosen
basis $\x$, these two maps are base-independent as
one can easily check (this will also follow from \ref{2ndDefOfKappa}).
Summarizing these results, we have the following theorem.
\end{Num}
\begin{Thm}
\label{ClassifyOrbits}
Let $\cL^{(3)}\subseteq\cL\times\cL\times\cL$ denote the set
of all triples of pairwise opposite Lagrangians. Then
we have a $\U(M)$-invariant surjective map
\[
\cL^{(3)}\rTo^\kappa \Herm_\eps^\circ(n)
\]
whose fibers are the $\U(M)$-orbits in $\cL^{(3)}$. The map $\kappa$ is given
by
\[
\kappa(g(X),g(Y),gu_t(Y))=\bra t,
\]
where $X,Y$ is our fixed pair of opposite Lagrangians as in
\ref{MatrixConvention}.
\qed
\end{Thm}
The result will be refined in Proposition \ref{ClassifyActionOnTriple}.
\begin{Num}
\label{CasesOrbitTypes}
According to \ref{ClassesOfgroups}, we have the following cases.

\medskip
\noindent
\textbf{Symplectic groups:} The triples are classified
by isomorphism classes of nondegenerate symmetric matrices.

\medskip
\noindent
\textbf{Hyperbolic orthogonal groups:}
The triples are classified by isomorphism classes of nondegenerate skew
symmetric matrices. There is one such class if $n$ is
even, and $\cL^{(3)}=\emptyset$ if $n$ is odd.

\medskip
\noindent
\textbf{Standard hyperbolic unitary groups:}
We may assume that $\eps=1$ (so the form is skew hermitian),
and then 
the triples are classified by isomorphism classes of $n$-dimensional
nondegenerate hermitian forms.
\end{Num}

\section{Flag complexes of graphs}
We continue to assume that $M$ is a $-\eps$-hermitian hyperbolic
module.
Now we consider the simplicial complex whose $k$-simplices
are $k+1$-sets of pairwise opposite Lagrangians. It
will be convenient to do this in the general setting of
graphs, flag complexes and simplicial sets.
\begin{Num}\textbf{The opposition graph}\\
By a \emph{graph} $\Gamma=(V,E)$ we understand an undirected graph
without loops or multiple edges; $V$
is its set of \emph{vertices}, 
$E$ its set of \emph{edges}, and edges are unordered pairs of
vertices.
If $\{u,v\}$ is an edge, we call $u,v$ \emph{adjacent}.
For the hyperbolic module $M$, we put $V=\cL$ and
$\cO=\left\{\{X,Y\}\mid X,Y\in\cL\text{ and }M=X+Y\right\}$.
The resulting graph $\Gamma=(\cL,\cO)$ is called the
\emph{opposition graph} of $M$.
\end{Num}
\begin{Num}\textbf{Flag complexes}\\
The \emph{flag complex} $Fl(\Gamma)$ of a graph $\Gamma$
is the simplicial set whose
$k$-simplices are tuples $(x_0,\ldots,x_k)$ of vertices,
such that for all $0\leq\mu<\nu\leq k$ we have
either $x_\mu=x_\nu$ or $\{x_\mu,x_\nu\}\in E$.
We have the standard $\ZZ$-free chain complex
$C_*(Fl(\Gamma))$ with the usual boundary operator
\[\textstyle
\partial(x_0,\ldots,x_k)=\sum_\nu(-1)^\nu
(x_0,\ldots,\hat x_\nu,\ldots,x_k)
\]
and the resulting homology and cohomology groups.
\end{Num}
We will also use \emph{alternating chains}, which are defined as
follows \cite{ES}. Let $N_k$ denote the submodule of
$C_k(Fl(\Gamma))$ generated by all elements $(x_0,\ldots,x_k)$
with $x_\mu=x_\nu$ for some $\mu<\nu$, and all
elements of the form
$(x_0,\ldots,x_k)-sign(\pi)(x_{\pi_0},\ldots,x_{\pi_k})$, for
$\pi\in Sym(k+1)$. The alternating chain complex 
is defined as the quotient chain complex
\[
\tilde C_*(Fl(\Gamma))=C_*(Fl(\Gamma))/N_*.
\]
The natural projection
$C_*(Fl(\Gamma))\rTo\tilde C_*(Fl(\Gamma))$ is a chain equivalence,
i.e. induces an isomorphism in homology and cohomology, 
see \cite{ES} VI.6.
The coset of $(x_0,\ldots,x_k)$ is denoted
$\bra{x_0,\ldots,x_k}$, with the relations
$\bra{x_0,\ldots,x_k}=0$ if $x_\mu=x_\nu$ for some $\mu<\nu$, and
\[
\bra{x_0,\ldots,x_k}=sign(\pi)\bra{x_{\pi_0},\ldots,x_{\pi_k}}.
\]
\begin{Num}\textbf{Equivariant cohomology}\\
The unitary group $\U(M)$ acts in a natural way
on the opposition graph and its flag complex.
In general, when a group $G$ acts (from the left, say)
on a chain complex $C_*$,
then we may consider the \emph{equivariant homology}
of $C_*$, which is defined as
follows. If $P_*\rTo\ZZ$ is a projective resolution of $G$ over
$\ZZ$, then the equivariant homology $H^G_*(C_*)$ is defined as the
total homology of the double complex $P_*\otimes_G C_*$,
see \cite[Ch.~VII.5]{Bro}.
The two canonical filtrations on the double complex
yield two spectral sequences $'\!E$ and $''\!E$
converging to $H_*^G(C_*)$ and the first one has on its second page
\[
'\!E^2_{pq}=H_p(G;H_q(C_*)).
\]
If $C_*$ is acyclic (eg., if $C_*=\ZZ$ is concentrated in dimension $0$),
then $'\!E$ collapses on the second page, and there is a
natural isomorphism $H_*^G(C_*)\cong H_*(G)$.

Similar remarks hold for cohomology; here, one looks at the double
complex $\Hom_G(P_*,C^*)$. Note also that if
$c:C_*\rTo A$ is a $G$-invariant cochain (so $G$ acts trivially on the
coefficient module $A$) and if $\eta:P_1\rTo\ZZ$ is the
augmentation map, then $c$ may be viewed in a natural way as a cochain in
$\Hom_G(P_*,\Hom_\ZZ(C_*,A))\cong \Hom_\ZZ(P_*\otimes_GC_*,A)$ via 
\[
c(p\otimes z)=\eta(p)c(z).
\]
\end{Num}
It is well-known that for a complete
graph (i.e. for $E={V\choose 2}$) the simplicial set
$FL(\Gamma)$ is acyclic. The following concept is a 
weakening of (infinite) complete graphs.
\begin{Num}\textbf{The star property}\\
A (nonempty)
graph $\Gamma=(V,E)$ has the \emph{star property} if for 
every finite set $x_0,\ldots,x_k$ of vertices, there exists
a vertex $y$ which is adjacent to the
$x_\nu$, for $\nu=0,\ldots,k$.
\end{Num}
Note that we require that $y\neq x_0,\ldots,x_k$. A graph with
the star property is obviously infinite. Note also that
the opposition graph of a hyperbolic module has
by \ref{EnoughOpposites} the star property 
if $D^\eps$ is infinite.
\begin{Lem}
If $\Gamma$ has the star property, then $Fl(\Gamma)$ is acyclic.

\proof If $(x_0,\ldots,x_k)$ is a $k$-simplex in
$C_k(Fl(\Gamma))$ and $y$ if
is adjacent to $x_0,\ldots,x_k$ put
$y\#(x_0,\ldots,x_k)=(y,x_0,\ldots,x_k)$. 
Suppose that $c$ is a $k$-cycle, i.e. $c$ is a finite linear
combination of $k$-simplices and $\partial c=0$.
Let $y$ be a vertex adjacent to all vertices appearing in the
simplices of $c$. Then $\partial(y\#c)=c-y\#\partial c=c$,
so $c$ is a boundary.
\qed
\end{Lem}
The geometric realization $|Fl(\Gamma)|$ of the flag complex
of a graph with 
the star property is in fact contractible. To see this, it
suffices by Hurewicz' Theorem to show that
$\pi_1|Fl(\Gamma)|=0$, see \cite[Ch.~7.6.24 and 7.6.25]{Spa}.
But from the star property,
any simplicial path in $|Fl(\Gamma)|$ is contained in a
contractible subcomplex, and every path is homotopic to a simplicial
path \cite[3.6]{Spa}.
\begin{Prop} If $D^\eps$ is infinite, then the flag complex of
the opposition graph of a $-\eps$-hermitian hyperbolic module
is acyclic and its geometric realization is
contractible. Consequently, we have in equivariant homology
a natural isomorphism
\[
H^{\U(M)}_*(Fl(\Gamma))\rTo^\cong H_*(\U(M))
\]
induced by the constant map $Fl(\Gamma)\rTo Fl(\{pt\})$,
and similarly for cohomology.
\qed
\end{Prop}
Using this natural isomorphism, we often identify these two
(co)homology groups.

\section{The projectivity groupoid}
If $X$ and $Y$ are opposite Lagrangians in the hyperbolic
module $M$, then we have canonical isomorphisms
$X\oplus X^J\cong X\oplus Y\cong Y\oplus Y^J$,
such that the first isomorphism is the identity on $X$ and the
second isomorphism is the identity on $Y$. In this way, we
associate an isomorphism $X\oplus X^J\rTo Y\oplus Y^J$
to every \emph{oriented} edge $(X\rTo Y)$
of the opposition graph $\Gamma$.
\begin{Num}\textbf{The projectivity groupoid}\\
\label{ProjectivityGroupoid}
Recall that a groupoid is a small category where every
arrow is an isomorphism. The \emph{projectivity groupoid}
$\cG M$ of $M$ is defined as follows. The objects 
of $\cG M$ are
$2$-graded vector spaces $X_*$ with
$X_1=X$ and $X_{-1}=X^J$, where
$X\in\cL$ is a Lagrangian. To each oriented edge $(X\rTo Y)$
we associate an isomorphism $[Y;X]:X_*\rTo Y_{-*}$
of degree $-1$, the composite
\[
[Y;X]:X\oplus X^J\rTo^\cong X\oplus Y\rTo^\cong Y\oplus Y^J.
\]
These maps generate the morphisms of $\cG M$.
We note that each object $X_*$ in $\cG$ carries a natural structure of a
hyperbolic module with $-\eps$-hermitian form $h_X$,
and that the morphisms preserve this structure. Furthermore 
\[
[X;Y][Y;X]=\id_{X*},
\]
so a morphism along a simplicial path depends only on the
homotopy class of the path in $\Gamma$ (i.e.
we have a natural transformation from the
\emph{fundamental groupoid} $\pi_1\Gamma$ to $\cG M$).
Finally, we note that $\cG M$ is in a natural way
$2$-graded: the paths of even length induce maps
of degree $1$, and the paths of odd length maps of
degree $-1$.
\end{Num}
\begin{Num}
\label{DetermineMorphisms}
Now we determine the morphism corresponding to a
closed path of length $3$. Let $X,Y$ be opposite
Lagrangians with bases $\x, \y$ as in \ref{MatrixConvention}
and let $Z=u_t(Y)$. We write $[Z;Y;X]=[Z;Y][Y;X]$ and so on.
Then 
\[
[Z;X;Y](\y_\nu)=u_t(\y_\nu)\quad\text{ and }\quad
[Z;X;Y](\x_\nu)=u_t(\x_\nu)=\x_\nu.
\]
Now
\[\textstyle
h(\y_\lambda,u_t(\y_\nu)))=h\left(\y_\lambda,\y_\nu+\sum_\mu\x_\mu
t_{\mu,\nu}\right)=
h\left(\y_\lambda,\sum_\mu\x_\mu t_{\mu,\nu}\right).
\]
The dual basis of $\y$ is $h(-,\x)^J$. With respect to the graded basis
$(\y,h(-,\x)^J)$
for $Y_*=Y\oplus Y^J$, the morphism $\phi=[Y;Z;X;Y]$ is therefore
given by a block matrix of the form
$\phi=\scriptstyle
\begin{pmatrix} *&*\\ t&*\end{pmatrix}$.
As this matrix has to be unitary and of degree $-1$, and because
$t^J=t\eps$, we obtain
\[
\phi=\scriptstyle
\begin{pmatrix} 0&-t^{-1}\\ t&0\end{pmatrix}.
\]
\end{Num}
\begin{Num}
\label{2ndDefOfKappa}
If $h_Y$ denotes the canonical $-\eps$-hermitian form on $Y_*$, then
\[
h_\phi(-,-)=h_Y(-,\phi(-))(-\eps)
\]
is the $\eps$-hermitian form
$h_\phi=
\scriptstyle
\begin{pmatrix} t&0\\0&t^{-J}\end{pmatrix}$. We note 
that $t^Jt^{-J}t=t$, so both blocks represent
the same isomorphism type $\bra t$ in $Herm_\eps^\circ(n)$,
and we define 
\[
\tilde\kappa (Z,X,Y)=\bra{t}.
\]
Note also that this class does not depend on the
basis $\y$ and that $\tilde\kappa$ is $\U(M)$-invariant. Furthermore,
we have $\tilde\kappa(Z,X,Y)=\kappa(X,Y,Z)$, where $\kappa$ is the 
invariant from Theorem \ref{ClassifyOrbits}.
We will see shortly that both invariants agree completely.
\end{Num}
\begin{Num}
\label{PropertiesOfKappa}
From $\phi^{-1}=\scriptstyle
\begin{pmatrix} 0&t^{-1}\\-t&0\end{pmatrix}$
we see that
\[
\tilde\kappa(X,Z,Y)=\bra{-t}.
\]
Next we note that for $y_1,y_2\in Y_*$ we have
\begin{eqnarray*}
h_Y(y_1,\phi(y_2))&=&h_X([X;Y]y_1,[X;Y]\phi[Y;X][X;Y]y_2)\\
&=&h_X([X;Y]y_1,[X;Y][Y,Z,X,Y][Y;X][X;Y]y_2)\\
&=&h_X([X;Y]y_1,[X;Y,Z,X][X;Y]y_2),
\end{eqnarray*}
whence
\[
\tilde\kappa(Y,Z,X)=\tilde\kappa(Z,X,Y),
\]
i.e. $\tilde\kappa$ is invariant under cyclic permutations
of the arguments. In particular, we have
\[
\tilde\kappa=\kappa
\]
Since $\kappa$ classifies by Theorem \ref{ClassifyOrbits}
triples of pairwise opposite Lagrangians,
we have the following sharpening of \ref{ClassifyOrbits}.
\end{Num}
\begin{Prop}
\label{ClassifyActionOnTriple}
The set-wise $\U(M)$-stabilizer of a triple $X,Y,Z$ of pairwise
opposite
Lagrangians induces (at least) the cyclic group $\ZZ/3$ on this
set. It induces the full symmetric group $Sym(3)$ if and only if
$a^Jta=-t$ for some $a\in\GL_nD$, where $\kappa(X,Y,Z)=\bra t$.
\qed
\end{Prop}

\section{The Maslov cocycle}
We want to turn the invariant
$\kappa:\cL^{(3)}\rTo\Herm_\eps^\circ(n)$ into
a $2$-cocycle for the flag complex $Fl(\Gamma)$ of the opposition graph.
Suppose that $A$ is an abelian group and that
$\alpha:\Herm_\eps^\circ(n)\rTo A$ is a map. By the properties
of $\kappa$ derived in \ref{PropertiesOfKappa} we see that
\[
c:\bra{X,Y,Z}\mapsto \alpha(\kappa(X,Y,Z))
\]
is a $2$-cochain
on the alternating chain complex $\tilde C_2(Fl(\Gamma))$,
provided that we have the relation $\alpha(\bra{-t})=-\alpha(\bra t)$
for all $t\in\Herm_\eps^\circ(n)$.
Now we investigate under what conditions this map is a cocycle, i.e.
under what conditions 
$c(\partial\bra{X,Y,Z,Z'})=0$, i.e. when
\[
c(\bra{Y,Z,Z'}-\bra{X,Z,Z'}+\bra{X,Y,Z'}-\bra{X,Y,Z})=0.
\]
\begin{Num}
We fix again $(X,\mathbf x),(Y,\mathbf y)$ as in \ref{MatrixConvention}.
Suppose that $Z=u_t(Y)$ and $Z'=u_{t'}(Y)$, and that
$X,Y,Z,Z'$ are pairwise opposite. 
So we have 
\[
\kappa(Z,X,Y)=\bra{t}\quad\text{ and }\quad
\kappa(Z',X,Y)=\bra{t'}.
\]
As $u^{-1}_t(Z)=Y$ and $u^{-1}_tu_{t'}=u_{-t+t'}$
we obtain
\[
\kappa(Z',X,Z)=\kappa(u^{-1}_tu_{t'}(Y),X,Y)=\bra{t'-t}.
\]
It remains to determine $\kappa(Z',Y,Z)$.
Let $w=\scriptstyle\begin{pmatrix} 0&1\\ -\eps&0\end{pmatrix}$.
Then $w$ is unitary and interchanges $X$ and $Y$.
We have $w(Z)=wu_t(Y)=wu_tw^{-1}(X)$ and we put
$v_t=wu_tw^{-1}=\scriptstyle
\begin{pmatrix} 1&0\\ -t\eps&1\end{pmatrix}$.
Then
\[
u_rv_t=
\scriptstyle
\begin{pmatrix} 1-rt\eps&r\\ -t\eps&1\end{pmatrix},
\]
whence $u_rw(Z)=u_rv_t(X)=Y$ for $r=t^{-1}\eps$.
So far we have achieved
\[
u_rw(Y)=X\quad\text{ and }\quad u_rw(Z)=Y.
\]
We seek $t''$ such that
$u_{t''}(Y)=u_rw(Z')=u_rwu_{t'}(Y)$, or $Y=u_{r-t''}wu_{t'}(Y)$.
Now
\[
u_{r-t''}wu_{t}'=\scriptstyle
\begin{pmatrix} 1&r-t''\\0&1\end{pmatrix}
\begin{pmatrix} 0&1\\ -\eps&-t'\eps\end{pmatrix}
=
\begin{pmatrix} (t''-r)\eps &1+(t''-r)t'\eps\\ -\eps&-t'\eps\end{pmatrix}
\]
whence $1=(r-t'')t'\eps$, which gives
$t''=r-{t'}^{-1}\eps=(t^{-1}-t'^{-1})\eps$, so
\[
\kappa(Z',Y,Z)=h(u_{t''}(Y),X,Y)=\bra{t^{-J}-t'^{-J}}.
\]
Plugging this into the boundary formula, we have the next result.
\end{Num}
\begin{Prop}
\label{CocycleCondition}
Let $A$ be an abelian group.
A function $\alpha:\Herm_\eps^\circ(n)\rTo A$ determines a
$\U(M)$-invariant
$2$-cocycle $c$ on the alternating $2$-chains 
of $Fl(\Gamma)$ if and only if
the following two relations hold for all $r,s,t\in \Herm_\eps^\circ(n)$:
\begin{alignat*}{2}
r+s=0  &\quad\text{ implies }\quad&\alpha\bra r+\alpha\bra{s}=0&\\
r+s+t=0&\quad\text{ implies }\quad&\alpha\bra r+\alpha\bra s+ \alpha\bra t
+\bra{-r^{-J}-s^{-J}}=0&
\end{alignat*}
\qed
\end{Prop}
Recall that the Grothendieck-Witt group $KU_0^\eps(D,J)$ of
hermitian forms is defined as the abelian group completion of the
commutative monoid consisting of the isomorphism classes of
nondegenerate
$\eps$-hermitian forms \cite[p.~239]{Sch}. The \emph{Witt group}
$W^\eps(D,J)$ is the  factor group of $KU_0^\eps(D,J)$ by the
subgroup generated by the $\eps$-hermitian hyperbolic modules.
We let $[t]$ denote the image of $\bra t$ in 
$KU_0^\eps(D,J)$ and $W^\eps(D,J)$.
\begin{Thm}
Let $\alpha\bra t=[t]\in W^\eps(D,J)$. Then $\alpha$ satisfies the
two conditions of Proposition \ref{CocycleCondition} and therefore
\[
\mas:\bra{X,Y,Z}\mapstoo [\kappa(X,Y,Z)]
\]
defines a
$W^\eps(D,J)$-valued $\U(M)$-invariant
$2$-cocycle on the alternating chain complex\linebreak
$\tilde C_2(Fl(\Gamma))$.

\proof
We proceed similarly as \cite[Prop 1.2]{PPS} and
use the fact that metabolic forms vanish in the Witt group
$W^\eps(D,J)$, see \cite[7.3.7]{Sch}, and that a $2k$-dimensional
nondegenerate hermitian form is metabolic if it admits a totally
isotropic subspace of dimension $k$.

For the $2n$-dimensional $\eps$-hermitian
form $(r)\oplus(-r)$, the vectors 
$(x,x)$, with $x\in D^n$, span an
$n$-dimensional totally isotropic subspace, so this form
is metabolic and $[r]+[-r]=0$.

Similarly we find for $r+s+t=0$ and the $4n$-dimensional form
$(r)\oplus(s)\oplus(t)\oplus(-r^J-s^J)$
that the vectors $(x,x,x,0)$ and $(r^{-J}x,s^{-J}x,0,x)$,
with $x\in D^n$, span a
totally isotropic $2n$-dimensional subspace, so this form is also
metabolic and 
$[r]+[s]+[t]+[-r^J-s^J]=0$.
\qed
\end{Thm}
\begin{Num}\textbf{The Maslov cocycle}\\
We call the $W^\eps(D,J)$-valued cocycle 
\[
\mas:\bra{X,Y,Z}\mapstoo [\kappa(X,Y,Z)]
\]
(and the corresponding cocycle for
the equivariant homology of $Fl(\Gamma)$) the \emph{Maslov cocycle}.
\end{Num}

\section{Naturality of the Maslov cocycle}
\label{Naturality}
We now study naturality of the Maslov cocycle under
restriction maps. There are two obvious types,
coming from field and from vector space inclusions.
We start with field inclusions, which are easier.
\begin{Num}\textbf{Extension of scalars}\\
Suppose that $D$ and $E$ are division rings
with involutions $J$ and $K$, respectively,
and that $\phi:D\rTo E$ is a homomorphism
commuting with these involutions.
If $M$ is a hyperbolic module over $D$, then
$M\otimes_\phi E$ is hyperbolic over $E$. The map
sending a Lagrangian $X\subseteq M$ to $X\otimes_\phi E$
induces an injection $\cL(M)\rTo \cL(M\otimes_\phi E)$
and an injection $\Gamma(M)\rTo\Gamma(M\otimes_\phi E)$
on the respective opposition graphs.
There is a natural map $W^D_E:W^\eps(D,J)\rTo W^\eps(E,K)$
and obviously, this map takes the Maslov cocycle $\mas_D$ of
$M$ to the Maslov cocycle $\mas_E$ of $M\otimes_\phi E$,
\begin{diagram}
\tilde C_2Fl(\Gamma(M)) & \rTo & \tilde C_2Fl(\Gamma(M\otimes_\phi E)\\
\dTo_{\mas_D} &&\dTo_{\mas_E}\\
W^\eps(D,J)&\rTo^{W^D_E}& W^\eps(E,K).
\end{diagram}
This gives the following result.
\end{Num}
\begin{Thm}
Let $\phi:(D,J)\rTo(E,K)$ be a homomorphism of skew fields
with involutions and assume that $D^\eps$ is infinite.
Consider the natural group monomorphism
\[
\Phi:\U(M)\rTo\U(M\otimes_\phi E).
\]
Then 
$(W^D_E)_*\mas_D=\Phi^*\mas_E$ in the diagram
\begin{diagram}
H^2(\U(M);W^\eps(D,J))&\rTo^{(W^D_E)_*} &
H^2(\U(M);W^\eps(E,K))\\
&&\uTo^{\Phi^*}\\
&& 
H^2(\U(M\otimes_\phi E);W^\eps(E,K)).
\end{diagram}
\qed
\end{Thm}
\begin{Num}
Suppose now that $M_1$ and $M_2$ are hyperbolic modules
(both over $D$) 
with corresponding sets $\cL_1,\cL_2$ of Lagrangians.
Then their direct sum $M=M_1\oplus M_2$ is in a natural way a
hyperbolic module. There is an obvious map
\[
\U(M_1)\rTo\U(M)
\]
and the question is what happens with the Maslov cocycle
under this map.
The problem is that the opposition graph
$\Gamma_1$ of $M_1$ is not a subgraph of the opposition graph
$\Gamma$ of $M$. However, there is a natural subgraph of
$\Gamma$ which projects $\U(M_1)$-equivariantly
onto $\Gamma_1$ and which yields a good comparison map.
The construction is as follows.

If $X_1\subseteq M_1$ and $X_2\subseteq M_2$ are Lagrangians,
then $X_1\oplus X_2$ is Lagrangian in $M$, so we have a
natural injection $\cL_1\times\cL_2\rTo\cL$.
Moreover, $X_1\oplus X_2$ is opposite $Y_1\oplus Y_2$ in
$M$ if and only if $X_\nu$ is opposite $Y_\nu$, for
$\nu=1,2$. This leads us to the following notion.
\end{Num}
\begin{Def}
The \emph{categorical product} $\Gamma_1\times\Gamma_2$
of two graphs has $V_1\times V_2$ as its set of vertices
and $(x_1,x_1)$ and $(y_1,y_2)$ are adjacent if and only if
$\{x_1,y_1\}\in E_1$ and $\{x_2,y_2\}\in E_2$. There are
natural maps
$\Gamma_1\lTo\Gamma_1\times\Gamma_2\rTo\Gamma_2$ with the
usual universal properties.
\end{Def}
The next result is immediate.
\begin{Lem}
The categorical product of two graphs having the
star property has again the star property.
In particular, its
flag complex is acyclic. \qed
\end{Lem}
Note that the categorical product of the graph consisting
of one single edge with itself is not even connected; 
\begin{diagram}[size=1em,abut]
\bullet&&&&&&&&\bullet&&\bullet \\
\dLine&\times&\bullet&\rLine&\bullet&& = &&&\ldLine\rdLine\\
\bullet&&&&&&&&\bullet&&\bullet\rlap{.}
\end{diagram}
The fact that $y\neq x_0,\ldots,x_k$ in the star property is crucial
for the Lemma.
\begin{Num}
So far we have for $\nu=1,2$ a diagram of $\U(M_1)$-equivariant
maps
\begin{diagram}[size=2em]
Fl(\Gamma_1)\\
& \luTo^{pr_1}\\
&&Fl(\Gamma_1\times\Gamma_2)&\rTo &Fl(\Gamma)\\
& \ldTo_{pr_2}\\
Fl(\Gamma_2)
\end{diagram}
and if $D^\eps$ is infinite, these three complexes are
acyclic.
Next we note that if we have a triangle 
$(X_1\oplus X_2,Y_1\oplus Y_2,Z_1\oplus Z_1)$
in $\Gamma_1\times\Gamma_2$ and if we choose bases
$\x_1,\x_2,\y_1,\y_2$ for $X_1,X_2,Y_1,Y_2$, then
\[
[\kappa(X_1\oplus X_2,Y_1\oplus Y_2,Z_1\oplus Z_1)]
=[t_1\oplus t_2]=[t_1]+[t_2],
\]
with $\kappa(X_\nu,Y_\nu,Z_\nu)=\bra{t_\nu}$.
Thus we have a commutative diagram
\begin{diagram}
\tilde C_2Fl(\Gamma_1) & \rTo^{i_1} & \tilde C_2Fl(\Gamma_1)\oplus 
\tilde C_2Fl(\Gamma_2)
&\lTo^{(pr_1,pr_2)} & \tilde C_2Fl(\Gamma_2\times\Gamma_2) \\
&\rdTo_{\mas_1}&\dTo_{\mas_1+\mas_2}&\ldTo_\mas\\
&&W^\eps(D,J).
\end{diagram}
which yields in cohomology
\begin{diagram}
H^2_{\U(M_1)}(Fl(\Gamma_1))&\lTo&
H^2_{\U(M_1)}(Fl(\Gamma_1))\oplus
H^2_{\U(M_1)}(Fl(\Gamma_2))&\rTo&
H^2_{\U(M_1)}(Fl(\Gamma_1\times\Gamma_2))\\
[\mas_1] &\lMapsto&[\mas_1]+[\mas_2]&\rMapsto& [\mas]
\end{diagram}
(we omit here the coefficient group $W^\eps(D,J)$).
Note that $[\mas_2]=0$ in $H^2_{\U(M_1)}(Fl(\Gamma_2))$,
as $\U(M_1)$ acts trivially on $Fl(\Gamma_2)$.
Mapping to the one-point space $\{pt\}$, we see that
$[\mas_1]$ and $[\mas]$ have the same image in
$H^2_{\U(M_1)}(\{pt\})=H^2(\U(M_1))$, and from
\begin{diagram}
H^2_{\U(M_1)}(Fl(\Gamma_1))&&
H^2_{\U(M_1)}(Fl(\Gamma_1\times\Gamma_2))&\lTo^\cong&
H^2_{\U(M_1)}(Fl(\Gamma))&\lTo&
H^2_{\U(M)}(Fl(\Gamma))\\
\uTo^\cong&&\uTo^\cong&&\uTo^\cong&&\uTo^\cong\\
H^2(\U(M_1))&\rEq&
H^2(\U(M_1))&\rEq&
H^2(\U(M_1))&\lTo&
H^2(\U(M)).
\end{diagram}
we obtain the following stability result.
\end{Num}
\begin{Thm}
Assume that $D^\eps$ is infinite,
let $M_1,M_2$ be hyperbolic modules and put $M=M_1\oplus M_2$.
Then the restriction map 
\[
H^2(\U(M_1);W^\eps(D,J))\lTo
H^2(\U(M);W^\eps(D,J))
\]
maps the Maslov cocycle $[\mas]$ for $\U(M)$
onto the Maslov cocycle $[\mas_1]$ for $\U(M_1)$.
\qed
\end{Thm}

\section{Reduction of the cocycle}

Our next aim is to show that the Maslov cocycle can be
reduced to a subgroup of the Witt group. 
For this, we need a refinement of the Lagrangians and
the opposition graph. We noted in \ref{CasesOrbitTypes}
that the Maslov
cocycle is trivial in the hyperbolic orthogonal situation,
where $J=\id$ and $\eps=-1\neq 1$, so we may
disregard this case.
By \ref{ClassesOfgroups} there is no loss of generality in
assuming that
\[
\eps=1
\]
in the remaining cases, and we will do this in this section.
\begin{Num}\textbf{Based Lagrangians}\\
Let $\Gamma=(V,E)$ be a graph and $f:X\rTo V$ a map. The
\emph{induced graph} $f^*\Gamma$
on $X$ is the graph whose vertices are the elements
of $X$, and $\{x,x'\}$ is an edge if and only if
$\{f(x),f(x')\}$ is an
edge of $\Gamma$. If $f$ is surjective and
if $\Gamma$ has the star property, then $f^*\Gamma$ also has the
star property. 
In what follows, we consider the set
$\widehat\cL$ of \emph{based Lagrangians}, i.e. pairs 
$(X,\x)$ where $X\subseteq M$ is a Lagrangian and $\x$
is a basis for $X$.
There is a forgetful surjection $F:\widehat\cL\rTo\cL$ and
we let 
\[
\widehat\Gamma=F^*\Gamma
\]
denote the induced graph on
this vertex set. We call $\widehat\Gamma$ the
\emph{based opposition graph}. Because the $\U(M)$-stabilizer
$P$ induces the full group $\GL(X)$ on $X$, we see that
$\U(M)$ acts transitively on $\widehat\cL$. With the notation
of \ref{MatrixConvention}, the stabilizer
of $(X,\x)$ is the group $U$. The map
$\widehat\Gamma\rTo\Gamma$ is equivariant, and
$Fl(\widehat\Gamma)$ is acyclic if $D^1$ is infinite.
In particular, we may use $Fl(\widehat\Gamma)$ to compute the
group cohomology of $\U(M)$.

We also have a
based version $\widehat\cG M$ of the projectivity groupoid.
The objects are again the $2$-graded spaces $X\oplus X^J$,
but now with a preferred graded basis consisting of $\x$ and the
dual basis of $\x$. The morphisms in $\widehat\cG M$ are
thus given by unitary matrices.
\end{Num}
\begin{Num}
\label{BasedCalculations}
We re-calculate the Maslov cocycle in terms of the based spaces.
In \ref{DetermineMorphisms}
we saw that we have in terms of our standard basis
$\x,\y$ the matrices
\begin{diagram}[width=5em]
(Y_*,\y) & \rTo^{\scriptsize
\begin{pmatrix}0&-1\\1&0\end{pmatrix}}
&(X_*,\x) & \lTo^{\scriptsize
\begin{pmatrix}0&-1\\ 1&0\end{pmatrix}}
&(Z_*,u_t(\y)) & \rTo^{\scriptsize
\begin{pmatrix}0&-t^{-1}\\ t&0\end{pmatrix}}
&(Y_*,\y).
\end{diagram}
If we add base changes through matrices
$a,b,c\in \GL_nD$ for $X$, $Y$ and $Z$ and reverse
the middle arrow, we arrive at the diagram
\begin{diagram}
(Y_*,b\y) & \rTo^{\scriptsize
\begin{pmatrix}0&-ab^J\\a^{-J}b^{-1}&0\end{pmatrix}}
&(X_*,a\x) & \rTo^{\scriptsize
\begin{pmatrix}0&ca^J\\ -c^{-J}a^{-1}&0\end{pmatrix}}
&(Z_*,cu_t(\y)) & \rTo^{\scriptsize
\begin{pmatrix}0&-bt^{-1}c^J\\ 
b^{-J}tc^{-1}&0\end{pmatrix}}
&(Y_*,b\y)
\end{diagram}
(and $cu_t(\y)=u_{c^{-J}tc^{-1}}(c\y)$).
With respect to the basis $b\y$, we have
\[
[Y;Z;X;Y]=
\begin{pmatrix}0&-bt^{-1}b^J\\ 
b^{-J}tb^{-1}&0\end{pmatrix}.
\]
Using invariants of these matrices, we now construct a 
refined cocycle.
\end{Num}
\begin{Num}\textbf{Invariants of hermitian forms}\\
The dimension induces a natural homomorphism 
$\dim:KU_0^1(D,J)\rTo\ZZ$.
Since the dimension of any hyperbolic module is even, there
is an induced map $W^1(D,J)\rTo\{\pm1\}$ mapping the
class $[t]$ to $(-1)^{\dim(t)}$. We denote its kernel by
$I(D,J)$; its elements are represented by even dimensional
hermitian forms. In the quadratic case ($J=\id$ and $\eps=1\neq-1$),
$I\!D=I(D,\id)$ is called the \emph{fundamental ideal}
in the Witt ring $W\!D=W^1(D,\id)$ \cite[Ch.~II.1]{Lam}.

Recall that the determinant is a homomorphism from $\GL_nD$ to
$K_1(D)$, the abelianization of $D^*=\GL_1D$.
The involution $J$ induces an automorphism $J$ on $K_1(D)$.
We let $N$ denote the subgroup of $K_1(D)$ consisting of elements
of the form $x^Jx$ and put $S=K_1(D)/N$. Since
$\det(g^Jtg)=\det(g^Jg)\det(t)$,
we have a well-defined homomorphism $[t]\mapstoo \det(t)N$ from
$KU_0^1(D,J)$ to $S$. However, this map cannot be factored through
$W^1(D,J)$. Similarly as in \cite[Ch.~II.2]{Lam} we introduce therefore
the abelian group
\[
\widehat S=S\times\{\pm1\},
\]
endowed with the commutative group law
\[
(x,(-1)^m)+(y,(-1)^n)=(xy(-1)^{mn},(-1)^{m+n}),
\]
and we define the \emph{signed discriminant} as
\[
disc(t)=(\det(t)N(-1)^{n(n-1)/2},(-1)^n),
\]
where $n=\dim(t)$. This map vanishes on hyperbolic forms and
induces therefore a homomorphism $disc:W^1(D,J)\rTo \widehat S$.
We let $I\!I(D,J)\subseteq W^1(D,J)$ denote the subgroup
generated by all elements $[t]$, where $\dim(t)\in 4\ZZ$ and
$\det(t)=1$. Obviously, $I\!I(D,J)\subseteq\mathrm{ker}(disc)$.
\end{Num}
\begin{Lem}
The sequence
\[
0\rTo I\!I(D,J)\rTo W^1(D,J)\rTo\widehat S
\]
is exact.

\proof
Let $[t]$ be a form in the kernel of $disc$. Then $\dim(t)$
is even and we distinguish two cases. If $\dim(t)=4$,
then $\det(t)=x^Jx\in N$. Choose $g\in\GL_nD$ with $\det(g)=x^{-1}$,
then $\det(g^Jtg)=1$ and $[t]=[g^Jtg]\in I\!I(D,J)$.
For $\dim(t)=4\ell+2$ we have
$\det(t)=-x^Jx$ and we consider the $4\ell+4$-dimensional form
$t\oplus h$, for
$h=\scriptsize\begin{pmatrix}0&-1\\1&0\end{pmatrix}$.
Then $\det(t\oplus h)=x^Jx$.
By the previous remark, $[t\oplus h]=[t]\in I\!I(D,J)$.
\qed
\end{Lem}
In the quadratic case, $I\!I(D,J)$ is the square 
$I^2D$ of the fundamental
ideal \cite[Ch.~II,2.1]{Lam}
\begin{Num}
\label{Reduction}
We define an $\widehat S$-valued equivariant $1$-cochain $f$ on
$\tilde C_1Fl(\widehat\Gamma)$ by
\[
f\bra{(X,a\x),(Y,b\y)}=(\det(-ab^J)(-1)^{n(n-1)/2}N,(-1)^n)
\in\widehat S,
\]
where the notation is as in \ref{BasedCalculations}.
Note that this is indeed an alternating cochain:
\[
(\det(g)N,(-1)^n)+(\det(-g^J)N,(-1)^n)=
(\det(-gg^J)(-1)^{n^2}N,(-1)^{2n})=(N,1).
\]
Then $df=f\partial$ is an $\widehat S$-valued $2$-coboundary,
and 
\begin{align*}
df\bra{(Z,c u_t(\y)),(X,a\x)),(Y,b\y)}&=
f\bra{(X,a\x),(Y,b\y))}
-f\bra{(Z,cu_t\y),(Y,b\y))}\\
&\qquad +f\bra{Z,c u_t(\y))(X,a\x))}\\
&=(\det(-ab^J)(-1)^{n(n-1)/2}N,(-1)^n)\\
&\qquad +(\det(-bt^{-1}c^J)(-1)^{n(n-1)/2}N,(-1)^n)\\
&\qquad +(\det(ca^J)(-1)^{n(n-1)/2}N,(-1)^n)\\
&=(\det(aa^Jbb^Jcc^Jt)(-1)^{n(n-1)/2}N,(-1)^n)\\
&=(\det(t)(-1)^{n(n-1)/2}N,(-1)^n)
\end{align*}
whence
\[
disc_*\mas +df=0
\]
where $disc_*$ denotes the coefficient homomorphism
induced by $disc :W^1(D,J)\rTo\widehat S$. 
Consequently, the image of $\mas$ vanishes in
$H^2(\U(M);\widehat S)$.
\end{Num}
\begin{Num}
Recall that $\mathbf{EU}(M)\subseteq\U(M)$ is the 
invariant subgroup generated by
the Eichler transformations. This group is perfect if
$D^1$ is infinite \cite[6.3.15]{HOM} and consequently
$H^1(\mathbf{EU}(M);A)=\Hom(\mathbf{EU}(M),A)=0$
for any coefficient group $A$ with trivial $\mathbf{EU}(M)$-action.
We put $\widehat S_0=disc(W^1(D,J))\subseteq\widehat S$.
The long exact cohomology sequences for the coefficient maps
\begin{diagram}
0&\rTo&\widehat S_0&\rTo&\widehat S&\rTo& \widehat S/\widehat S_0&\rTo& 1\\
0&\rTo&I\!I(D,J)&\rTo& W^1(D,J)&\rTo&\widehat S_0&\rTo&0
\end{diagram}
yield therefore monomorphisms
\begin{diagram}
0&\rTo&H^2(\mathbf{EU}(M);\widehat S_0)&\rTo&
H^2(\mathbf{EU}(M);\widehat S)\\
0&\rTo&H^2(\mathbf{EU}(M);I\!I(D,J))&\rTo&
H^2(\mathbf{EU}(M);W^1(D,J)).
\end{diagram}
\end{Num}
This gives us the next Theorem. To keep notation simple,
we denote the restriction of $\mas$ to the subgroup
$\mathbf{EU}(M)$ also by $\mas$.
\begin{Thm}
Assume that $\eps=1$ and that $D^1$ is infinite.
There exists a unique cohomology class
$[\widetilde\mas]\in H^2(\mathbf{EU}(M);I\!I(D,J))$ which maps
under the coefficient homomorphism $I\!I(D,J)\rTo W^1(D,J)$
onto $[\mas]$. We call this class the \emph{reduced Maslov cocycle}.

\proof
As we proved in \ref{Reduction}, $disc_*[\mas]+[df]=0$ in
$H^2(\mathrm{EU}(M);\widehat S)$, whence
$disc_*[\mas]=0$ in $H^2(\mathbf{EU}(M);\widehat S_0)$.
Therefore
$[\mas]$ has a preimage $[\widetilde\mas]$ in
$H^2(\mathbf{EU}(M);I\!I(D,J))$. The map
$H^2(\mathbf{EU}(M);I\!I(D,J))\rTo H^2(\mathbf{EU}(M);W^1(D,J))$
is injective, so the preimage is unique.
\qed
\end{Thm}
\begin{Num}
\label{ReducedCocycleSymplectic}
In the symplectic situation  $(J,\eps)=(\id,1)$ it is possible
to give an explicit formula for the reduced cocycle
$\widetilde\mas$. Then $\widehat S=\widehat S_0$ and
$\mathbf{EU}(M)=\U(M)=\Sp_{2n}D$ and we can directly define a
$W^1(D,J)$-valued $1$-cochain on $Fl(\widehat\Gamma)$ by
\[
\tilde f\bra{(X,a\x),(Y,b\y)}=
\bra{\det(-ab),1,\ldots,1},
\]
where the right-hand side denotes as usual the $n$-dimensional
symmetric bilinear form with the given entries on the diagonal.
Under the map $p:Fl(\widetilde\Gamma)\rTo Fl(\Gamma)$ this
is a lift of $f$ and we have
$disc_*d\tilde f=p^*df$. Thus
\[
\widetilde\mas=p^*\mas+d\tilde f
\]
is the reduced Maslov
cocycle on $FL(\widehat\Gamma)$ in the symplectic case.
Explicitly, it reads as
\begin{align*}
\widetilde\mas\bra{(X,a\x),(Y,b\y),(Z,cu_t\y)}&=
\bra{\det(-ab),1,\ldots,1}+\bra{\det(ca),1,\ldots,1}\\
&\qquad+\bra{\det(-btc),1,\ldots,1}-\bra{t}.
\end{align*}
\end{Num}

\section{Kashiwara's Maslov cocycle}
\label{ClassicalMaslov}
In the symplectic situation over a field $D$ of characteristic
$\neq 2$, the Maslov index is classically
defined through a different quadratic form \cite{LV}.
(A variant is used in \cite{PPS}, while a topological generalization
for bounded symmetric domains of tube type is given in \cite{NO}.
See \cite{CLM} for a survey of topological Maslov indices.)
\begin{Num}\textbf{Kashiwara's Maslov index}\\
Let $D$ be a field of characteristic $\neq 2$. We assume that we are
in the symplectic situation $\eps=1$, $J=\id$.
Given three Lagrangian $X,Y,Z$ (not necessarily pairwise opposite)
we consider the following $3n$-dimensional quadratic
from $q_{X,Y,Z}$ on the direct sum $X\oplus Y\oplus Z$:
\[
q_{X,Y,Z}(x,y,z)=h(x,y)+h(y,z)+h(z,x).
\]
If the Lagrangians are not pairwise opposite, the quadratic form
is going to have a radical. The \emph{Kashiwara-Maslov index} of $(X,Y,Z)$
is the class in the Witt group $W\!D$ which is
represented by the nondegenerate part 
$q_{X,Y,Z}^+$ of $q_{X,Y,Z}$.

For $D=\RR$,
the Witt group $W\RR$ is isomorphic to $\ZZ$ via the signature
and the Maslov index can directly be defined as the signature of
$q_{X,Y,Z}$ (even if the form is degenerate). This is essentially
Kashiwara's definition of the symplectic
Maslov index as developed in \cite[1.5.1]{LV}.

If $X,Y,Z$ are pairwise opposite, we find that with respect to our
standard basis $\x,\y,u_t\y$ for $X\oplus Y\oplus Z$
the quadratic form is represented by the matrix
\[
q_{X,Y,Z}=
\begin{pmatrix}0&-1&0\\0&0&t\\
1&0&0\end{pmatrix}.
\]
We note that $X\oplus Y\oplus0$ is a hyperbolic submodule
in $X\oplus Y \oplus Z$
whose orthogonal complement is spanned by
vectors of the form $(tz,z,z)\in D^{3n}$. The restriction of
$q_{X,Y,Z}$ to this subspace is given by $z\mapstoo(z^Ttz)$,
so $q_{X,Y,Z}=q_{X,Y,Z}^+$ is represented by $[t]$ in $W\!D$.
This is our first result.
\end{Num}
\begin{Prop}
If $X,Y,Z$ are pairwise opposite Lagrangians, then
the Kashiwara-Maslov index of $(X,Y,Z)$ agrees with the image $[t]$
of $\bra{t}=\kappa(X,Y,Z)$ in the Witt group $W\!D$.
\qed
\end{Prop}
Next we get to Kashiwara's Maslov cocycle, which is defined 
as follows. We fix a Lagrangian $X_0\in\cL$ and define
$\tau:\Sp_{2n}D\times\Sp_{2n}D\rTo W\!D$ via
\[
\tau(g,h)=\bra{q_{X_0,g(X_0),gh(X_0)}^+}.
\]
We want to relate this group cocycle to our Maslov
cocycle defined in terms of the flag complex of the opposition
graph. 
\begin{Num}
Recall the bar notation \cite[I.5]{Bro} for the standard free resolution of a
group $G$ over $\ZZ$. Its chain complex is given as
\[
F_n=\ZZ G^{n+1}
\]
and the generator
$(1,g_1,g_1g_2,g_1g_2g_3,\ldots,g_1\cdots g_n)\otimes 1\in F_n\otimes_G\ZZ$
is denoted
$[g_1|\ldots|g_n]$. Then $\tau$ can be viewed as the $W\!D$-valued $2$-chochain
$[g|h]\mapsto \tau(g,h)$ for $G=\Sp_{2n}D$ and one verifies 
the cocycle identity \cite[1.5.8]{LV}   .
\end{Num}
In general, suppose that $X$ is a set on which a group $G$ acts,
and that $c:X\times X\times X\rTo A$ is a $G$-invariant map
taking values in an abelian group $A$,
such that $c$ satisfies the cocycle identity
$c(x,y,z)-c(w,y,z)+c(w,x,z)-c(w,x,y)=0$. If we choose
a base point $o\in X$, it is not difficult to
see that the cocycle
$(g_1,g_2,g_3)\mapstoo c(g_1(o),g_2(o),g_3(o))$
defined on the standard free resolution $F_*$ of $G$ over
$\ZZ$ and the cocycle
$g\otimes(x,y,z)\mapstoo c(x,y,z)$ defined on
$F_0\otimes_G C_2\subseteq F_*\otimes_G C_*$ are homologous
($C_*$ is the standard complex of $k+1$-tuples of elements of
$X$). However, we cannot use this directly to compare our
Maslov cocycle with its classical counterpart, since our
cocycle is defined only on special triples of Lagrangians.
We need to refine this idea, using some elementary homological
algebra. We do this in general, as we need it also in the
next section.
\begin{Num}
\label{StarResolution}
Let $\Gamma=(V,E)$ be a graph with the star property. Suppose
that $G$ is a group acting transitively on the vertices of $\Gamma$.
Let $o\in V$ be a base point and consider the induced graph
$\Gamma_G$ on $G$ under the map $G\rTo V$, $g\mapstoo g(o)$ and
its flag complex 
\[
F_*'=C_*Fl(\Gamma_G)\subseteq F_*.
\]
Obviously, this chain complex is a
free resolution of $G$ over $\ZZ$ and a subcomplex of the
standard free resolution $F_*$ of $G$.
Both chain complexes $F_*$ and $F_*'$ can be used to
determine the group (co)homology of $G$.

Suppose now that $c:C_2Fl(\Gamma)\rTo A$ is a $G$-invariant
cocycle. Then we can construct two $2$-cocycles for $G$, one
via 
\[
\hat c:(g_0,g_1,g_2)\mapstoo c(g_0(o),g_1(o),g_2(o))
\]
on  $F_2'\otimes_G\ZZ$, and the other via 
\[
c:g\otimes(x,y,z)\mapstoo c(x,y,z)
\]
on $F'_0\otimes_GC_2Fl(\Gamma)\subseteq F_*'\otimes_GC_2Fl(\Gamma)$.
Our first aim is to prove that both cocycles are
homologous.
We put $C_*=C_*Fl(\Gamma)$ and we call a generator
$(g_0,\ldots,g_m)\otimes(x_0,\ldots,x_n)\in F_m'\otimes_GC_n$
\emph{admissible} if $\{g_0(o),\ldots,g_m(o),x_0,\ldots,x_n\}$
consists of pairwise adjacent elements in $\Gamma$. This is
a well-defined notion, i.e. invariant under the left diagonal
action of $G$.
Let
$D_{**}\subseteq F_*'\otimes_G C_*$ denote the submodule
generated by the admissible elements. We note that this
submodule is $\ZZ$-free and closed under the vertical and horizontal
differentials, so it is a double complex.
\end{Num}
\begin{Lem}
\label{TrickyLemma}
The inclusion $D_{**}\rInto F_*'\otimes_G C_*$ induces
an isomorphism in homology and cohomology (for coefficient groups
with trivial $G$-action).

\proof
We show that the relative homology groups of the pair
$(F_*'\otimes_G C_*,D_{**})$ vanish. Let
$z\in\bigoplus_{m+n=k}F_m'\otimes_GC_n$ be a relative $k$-cycle and let
$\tilde z\in\bigoplus_{m+n=k}F_m'\otimes_\ZZ C_n$
be an element which maps onto $z$. We choose a group element $j$
such that for all terms
$(g_0,\ldots,g_m)\otimes (x_0,\ldots,x_n)$
appearing in $\tilde z$, the vertex $j(o)$ is adjacent to
$g_0(o),\ldots,g_m(o),x_0,\ldots,x_m$
(this is a well-defined condition as we work with
$\tilde z\in\bigoplus_{m+n=k}F_m'\otimes_\ZZ C_n$ where
the $G$-action is \emph{not} factored out).
Consider the $k+1$-chain $j\#\tilde z$,
whose $(m+1,n)$-terms are of the form
$(j,g_0,\ldots,g_m)\otimes (x_0,\ldots,x_n)$.
The total differential is
\[
\partial(j\#\tilde z)=\tilde z-j\#(\partial\tilde z).
\]
Projecting this equation back to 
$F_m'\otimes_GC_n$, we see that
the image of $j\#\partial\tilde z$ is in $D_{*+1,*}$.
Thus $z$ is a relative boundary and 
$H_*(F_*'\otimes_GC_*,D_{**})=0$.
From the long exact homology sequence we get an isomorphism
$H_*(D_{**})\rTo^\cong H_*(F_*'\otimes_GC_*)$.
Since both $F_*'\otimes_GC_*$ and
$D_{**}$ are $\ZZ$-free, the universal coefficient theorems and
the $5$-Lemma yield isomorphisms for homology and cohomology with
arbitrary coefficient groups $A$ (with trivial $G$-action),
see \cite[5.3.15,5.5.3]{Spa}.
\qed
\end{Lem}
The remaining part of the comparison is routine.
We denote elements of $G$ by $g,h,i$ and vertices of
$\Gamma$ by $u,v,w$.
We define two $1$-cochains $f_1,f_2$ on $D_{**}$
by 
\[
f_1((g)\otimes(u,v))=c(g(o),u,v)\quad\text{and}\quad
f_2((g,h)\otimes(u))=c(g(o),h(o),u),
\]
where $c$ is the given $G$-invariant $2$-cocycle on $Fl(\Gamma)$.
Then $df_\nu=f_\nu\partial$ and using the cocycle identity for $c$,
we obtain
\begin{align*}
df_1((g)\otimes(u,v,w))&=(c(g(o),v,w)-c(g(o),u,w)+c(g(o),u,v)\\
&=c(u,v,w)\\
df_1((g,h)\otimes(u,v))&=c(h(o),u,v)-c(g(o),u,v)\\
df_2((g,h)\otimes(u,v))&=-c(g(o),h(o),v)+c(g(o),h(o),u)\\
&=df_1((h,i)\otimes(x,y))      \\
df_2((g,h,i)\otimes(u))&=
c(h(o),i(o),u)-c(g(o),i(o),u)+c(g(o),h(o),u)\\
&=c(g(o),h(o),i(o))
\end{align*}
which shows that 
\[
df_1-df_2=c-\hat c.
\]
\begin{Thm}
Let $G$ be a group acting vertex-transitively on a graph
$\Gamma$ having the star property, let 
$c:Fl_2(\Gamma)\rTo A$ be a $G$-invariant $A$-valued $2$-cocycle
(where $G$ acts trivially on $A$).
Fix a vertex $o$ of $\Gamma$ and let
$F'_*\subseteq F_*$ and $C_*$ be as in \ref{StarResolution}.
Then the cocycles
\[
\hat c:F'_2\otimes_G\ZZ\rTo A,\qquad (g_0,g_1,g_2)\otimes 1
\mapstoo c(g_0(o),g_1(o),g_2(o))
\]
and
\[
c:F_0\otimes_G C_2\rTo A,\qquad g\otimes(x,y,z)\mapsto c(x,y,z)
\]
are homologous under the isomorphism
\[
H^2(G;A)\rTo^\cong H^2_G(C_*;A).
\]
Moreover, there exists a cocycle
$\hat{\hat c}:F_2\otimes_G\ZZ\rTo A$ extending $\hat c$, i.e.
$\hat c=\hat{\hat c}|_{F'_*}$.

\proof
Only the last claim remains to be proved. Since the
inclusion $F_*'\subseteq F_*$
induces an isomorphism in cohomology, we find a cocycle
$\tilde{\hat c}$ on $F_*\otimes_G\ZZ$ such that
$\hat c-\tilde{\hat c}|_{F'_*\otimes_G\ZZ}=da$ is
a coboundary. Now $F'_*\otimes_G\ZZ$ is a direct summand in the
$\ZZ$-free module
$F_*\otimes_G\ZZ$, so we can extend $a$ to a $1$-cochain $\tilde a$
on $F_*\otimes_G\ZZ$. Then
$(\tilde{\hat c}+d\tilde a)|_{F'_*\otimes_G\ZZ}=\hat c$.
\qed
\end{Thm}
\begin{Cor}
For a field $D$ of characteristic $\neq 2$,
Kashiwara's Maslov cocycle and our Maslov cocycle yield the
same cohomology class in $H^2(\Sp_{2n}D;W\!D)$.
\qed
\end{Cor}
We obtain also the following general result for unitary groups
over arbitrary skew fields.
\begin{Cor}
If $D^\eps$ is infinite and $o\in\cL$ is a fixed Lagrangian, then
there exists a group cocycle $\tau:\U(M)\times\U(M)\rTo W^1(D,J)$
such that 
\[
\tau(g,h)=\bra{\kappa(o,g(o),gh(o))}
\]
holds for all pairs $g,h$ with $o,g(o),gh(o)$ pairwise opposite.
\qed
\end{Cor}

\section{The Maslov cocycle as a central extension}

The reduced Maslov cocycle defines a central extension
\cite[IV.3]{Bro} \cite[1.4C]{HOM}
\[
1\rTo I\!I(D,J)\rTo \widehat{\mathbf{EU}(M)}\rTo\mathbf{EU}(M)\rTo1
\]
of $\mathbf{EU}(M)$ by $I\!I(D,J)$. This extension 
is uniquely determined by the homomorphism
\[
[\widetilde\mas]\in H^2(\mathbf{EU}(M);I\!I(D,J))\cong
\Hom(H_2(\mathbf{EU}(M)),I\!I(D,J));
\]
our aim is to determine this homomorphism
$H_2(\mathbf{EU}(M))\rTo I\!I(D,J)$ algebraically.
In view of the naturality we proved in Section \ref{Naturality},
we begin with the smallest case $\Sp_2D=\SL_2D$, where
$D$ is an infinite field. We do allow fields of characteristic
$2$, as we rely on results in \cite{Mat} and \cite{Moore}
which are valid over arbitrary (infinite) fields.
Note, however, that in our set-up
the Witt group $W^1(D,\id)$ is always the Witt group of
symmetric bilinear forms (and not of quadratic forms).

\begin{Num}\textbf{The Schur multiplier of $\SL_2D$ and the
Steinberg cocycle}\\
We put
\[
u_t=\begin{pmatrix}1&t\\0&1\end{pmatrix}
\qquad
a_r=\begin{pmatrix}r&0\\0&r^{-1}\end{pmatrix}
\qquad
b_r=\begin{pmatrix}0&r\\-r^{-1}&0\end{pmatrix}
\]
Since $\SL_2D$ is a two-transitive group, every element
is either of the form $a_ru_t$ or of the
form $u_sb_ru_t$.
We define $KSp_2D$ as the abelian group generated by
symbols  $\{x,y\}$, for
$x,y\in D^*$, (the \emph{symplectic Steinberg symbols}),
subject to the relations
\begin{align*}
\{st,r\}+\{s,t\}&=\{s,tr\}+\{t,r\}\\
\{s,1\}&=\{1,s\}=0\\
\{s,t\}&=\{t^{-1},s\}\\
\{s,t\}&=\{s,-st\}\\
\{s,t\}&=\{s,(1-s)t\}\quad\text{ if }s\neq 1.
\end{align*}
According to \cite[p.~199]{Moore} \cite[5.11]{Mat}
the Schur multiplier of
$\SL_2D$ is $H_2(\SL_2D)\cong KSp_2D$. Moreover,
the Steinberg normal form of the universal group cocycle 
\[
stbg:\SL_2D\times\SL_2D\rTo H_2(\SL_2D)
\]
is given for 'generic' group elements by
\[
stbg(g(s_1,r_1,t_1),g(s_2,r_2,t_2))=
\{{\textstyle\frac{t}{r_1r_2},-\frac{r_1}{r_2}}\}-\{-r_1,-r_2\}
% \{-r_1t^{-1},r_2\}-
% \{-r_1t^{-1},-t^{-1}\},
\]
where $t=t_1+s_2\neq 0$ and $g(s,r,t)=u_sb_ru_t$,
cp.~\cite[p.~198 (1)]{Moore}, \cite[5.12]{Mat}
and \cite{Kub} in a more special situation. We note that the formula given in
\cite[p.~198 (1)]{Moore} is incorrect. The formula above is due to Schwarze
\cite[5.9]{Schwarze} and agrees with Matsumoto's calculations.
\end{Num}
\begin{Num}
\label{ReductionMap}
Given $x,y\in D$, we denote by
$(x,y)_D$ the $4$-dimensional symmetric bilinear form
\[
(x,y)_D=\bra{1,-x,-y,xy}.
\]
If $\mathrm{char}(D)\neq 2$, this is the norm form of the
quaternion algebra $\left(\frac{x,y}D\right)$ \cite[2.\S11]{Sch}.
Obviously, $(x,y)_D\in I\!I(D,\id)$, and
$(x,y)_D=(y,x)_D=(xz^2,y)_D$. Using the fact that
the metabolic form $\bra{x,-x}$ vanishes in $W^1(D,\id)$, it is routine
to verify that these elements satisfy the first four defining
relations of $KSp_2D$; for example
$(s,-st)_D=\bra{1,-s,st,-s^2t}\cong\bra{1,-s,st,-t}\cong(s,t)_D$. 
For the last relation, it suffices
to check that $\bra{-t,st}\cong\bra{-(1-s)t,(1-s)st}$ for $s\neq 1$.
This follows from
\[\scriptsize
\begin{pmatrix}
1& 1\\s&1
\end{pmatrix}
\begin{pmatrix}
-t& 0\\0&st
\end{pmatrix}
\begin{pmatrix}
1& s\\1&1
\end{pmatrix}=
\begin{pmatrix}
-(1-s)t& 0\\0&(1-s)st
\end{pmatrix}
\]
Thus we have a homomorphism
\[
R:KSp_2(D)\rTo I\!I(D,\id)\subseteq W^1(D,\id)
\]
which maps
the symplectic Steinberg symbol $\{u,v\}$ to the $4$-dimensional
symmetric bilinear form $R(\{u,v\})=(u,v)_D$.
\end{Num}
Applying $R$ to the Steinberg cocycle, we obtain
(with the same notation as before) for 'generic' group elements
\begin{align*}
R\circ stbg(g(s_1,r_1,t_1),g(s_2,r_2,t_2))
&=({\textstyle\frac{t}{r_1r_2},-\frac{r_1}{r_2}})_D-
  (-r_1,-r_2)_D\\
&=(r_1r_2t,r_1r_2)_D-(-r_1,-r_2)_D\\
&=(-r_1r_2,t)_D-(-r_1,-r_2)_D\\
&=\bra{1,r_1r_2,-t,-r_1r_2t}-\bra{1,r_1,r_2,r_1r_2}\\
&=\bra{1,r_1r_2,-t,-r_1r_2,t,-1,-r_1,-r_2,-r_1r_2}\\
&=\bra{-t,-r_1r_2t,-r_1,-r_2}=-\bra{t,r_1r_2t,r_1,r_2}.
\end{align*}
% \begin{align*}
% R\circ stbg(g(s_1,r_1,t_1),g(s_2,r_2,t_2))
% &=(-r_1t^{-1},r_2)_D-
%   (-r_1t^{-1},-t^{-1})_D\\
% &=(-r_1t,r_2)_D-(-r_1t,-t)_D\\
% &=(-r_1t,r_2)_D-(-r_1,-t)_D\\
% &=\bra{1,r_1t,-r_2,-r_1tr_2}-\bra{1,r_1,t,r_1t}\\
% &=\bra{1,r_1t,-r_2,-r_1tr_2,-1,-r_1t,-t,-r_1t}\\
% &=\bra{-r_2,-r_1r_2t,-r_1,-t}\\
% &=-\bra{r_2,r_1r_2t,r_1,t}.
% \end{align*}
\begin{Num}
\label{Comparison}
We compare this expression with the reduced Maslov cocycle.
In $M=D^2$ we put $\x={1\choose0}$, $X=\x D$, and $o=(X,\x)$. 
Using the notation of \ref{StarResolution}, we have
for $F_2'$ the formula
\[
\tau(g_1,g_2)=\tau([g_1|g_2])=\widetilde\mas(o,g_1(o),g_1g_2(o))
\]
where three vertices $o,g_1(o),g_1g_2(o)$ have to be pairwise opposite.
For the first pair of vertices, this condition
gives  $g_1=u_{s_1}b_{r_1}u_{t_1}$, and for the
second pair  $g_2=u_{s_2}b_{r_2}u_{t_2}$. Then
\begin{align*}
\tau([g_1|g_2])
&=\widetilde\mas(o,g_1(o),g_1g_2(o))\\
&=\widetilde\mas(g_1^{-1}(o),o,g_2(o))\\
&=-\widetilde\mas(o,g_1^{-1}(o),g_2(o))\\
&=-\widetilde\mas(o,g(-t_1,-r_2,-s_1)(o),g(s_2,r_2,t_2)(o))\\
&=-\widetilde\mas(o,u_{-t_1}b_{-r_1}(o),u_{s_2}b_{r_2}(o))\\
&=-\widetilde\mas(o,b_{-r_1}(o),u_{t_1+s_2}b_{r_2}(o))
\end{align*}
which yields the additional condition $t=t_1+s_2\neq0$
that ensures that the first and third vertex are opposite.
Note that by \ref{StarResolution} the class of any $2$-cocycle
is completely determined by its values on 
$F_2'$, so it suffices indeed to work with 'generic' elements.
The explicit formula in \ref{ReducedCocycleSymplectic} for the
reduced Maslov cocycle yields now
$a=1$, $b=r_1^{-1}$ and $c=-r_2^{-1}$, whence
\[
\tau([g_1|g_2])=-\bra{t,r_1,r_2,r_1r_2t}=R_*stbg([g_1|g_2]).
\]
\end{Num}
For $\SL_2D$ over fields of characteristic $\neq 2$,
the following result was proved in \cite[Sec.~5]{Ne} and \cite{Ba}.
\begin{Thm}
\label{CocycleTheorem}
Let $D$ be an infinite field. The central extension of $\Sp_{2n}D$ 
determined
by the reduced Maslov cocycle is given by the homomorphism
$R:KSp_2D\rTo I\!I(D,\id)$.

\proof
For $n=1$ we showed this in \ref{Comparison} above.
In general, the standard inclusion $\Sp_{2n}D\rInto\Sp_{2n+2}D$
induces for all $n\geq 1$ an isomorphism in
$2$-dimensional homology, such that the
universal Steinberg cocycle for $\Sp_{2n+2}D$ restricts to the
universal Steinberg cocycle for $\Sp_{2n}D$
\cite[5.11]{Mat}. The result follows now by induction on $n$.
\qed
\end{Thm}
For fields of characteristic $\neq 2$, this is stated 
in \cite[3.1]{PPS}. However, the proof has a gap:
the authors evaluate the reduced Maslov cocycle on the torus
(the diagonal matrices)
and compare it there with the universal Steinberg cocycle.
But they fail to show that the reduced Maslov cocycle is
a Steinberg cocycle, so they cannot use the comparison theorem
\cite[5.10]{Mat}. 

In any case, this result settles the situation for symplectic 
groups over infinite fields of arbitrary characteristic.
Note that for fields of characteristic $\neq2$, the map $R$ is
surjective \cite[4.5.5]{Sch}, so $\widehat{\Sp_{2n}D}$ is an
epimorphic image of the universal central extension.
\begin{Num}\textbf{Local fields}\\
By a local field we mean a locally compact (nondiscrete) field;
the connected local fields are $\RR,\CC$ and the totally
disconnected ones are the finite extensions of the
$p$-adic fields $\mathbb Q_p$ and, in positive characteristic, the
fields $\mathbb F_q((X))$ of formal
Laurent series over finite fields \cite[1.3]{Weil}.
Being a closed subgroup of the general linear group,
a symplectic or unitary group over a local field is in a natural
way a locally compact group.
\end{Num}
\begin{Num}\textbf{The Maslov cocycle over $\RR$}\label{RealMaslov}\\
For $D=\RR$, the Witt group $W\RR=W^1(\RR,\id)$
is isomorphic to $\ZZ$ via the
signature $sig:W\RR\rTo\ZZ$ \cite[2.4.8]{Sch}; the fundamental
ideal $I\RR$ has index $2$, and
$I\!I\RR=I^2\RR$ has index $4$. We note that
\[
sig((x,y)_D)=\begin{cases} 4 &\text{ if }x,y<0\\
0&\text{else}\end{cases}
\]
By \cite[10.4]{Moore} \cite[p.~51]{Mat}, this $4\ZZ$-valued
cocycle yields
precisely the universal covering group 
$\widetilde{\Sp_{2n}\RR}$ of $\Sp_{2n}\RR$.

We compare the relevant classifying spaces.
Let $B\Sp_{2n}\RR^\delta$ denote the classifying space for
$\Sp_{2n}\RR$, viewed as a discrete topological group,
and $B\Sp_{2n}\RR$ the classifying space for the Lie group
$\Sp_{2n}\RR$; the latter is homotopy equivalent to
$B\U(n)$, as $\U(n)\subseteq\Sp_{2n}\RR$ is by
\cite[X~Tab.~V]{He} and Iwasawa's Theorem 
\cite[VI~\S2]{He} a homotopy equivalence.
The classifying space $B\Sp_{2n}\RR^\delta$
is an Eilenberg-MacLane space of type
$K(\Sp_{2n}\RR,1)$ whose cohomology is naturally isomorphic
to the abstract group cohomology of $\Sp_{2n}\RR$ \cite[II.4]{Bro}.
The identity map from the discrete group to the Lie group
induces a continuous map between the classifying spaces
\[
F:B\Sp_{2n}\RR^\delta\rTo B\Sp_{2n}\RR.
\]
On the right, the
universal covering is classified by the first Chern class $c_1$.
This shows that under the
forgetful map $F^*$ the first universal Chern class
$c_1\in H^2(B\SU(n))\cong H^2(B\Sp_{2n}\RR)$ pulls back to the Maslov cocycle,
\[
F^*(c_1)=[\widetilde\mas]
\]
(if the sign for the Chern classes is chosen appropriately).
The real
Maslov cocycle may be viewed therefore as a combinatorial
description of the first Chern class; this was observed in \cite{Tu}.
\end{Num}
\begin{Prop}
Under the forgetful functor from topological groups to abstract
groups, the first Chern class for $\Sp_{2n}\RR$ maps to the
reduced Maslov cocycle.
\qed
\end{Prop}
As $I\!I(\CC,\id)=0$ the Maslov cocycle for $\Sp_{2n}\CC$
vanishes. Now we turn to nonarchimedean local fields,
cp.~\cite[p.~104-115]{LV}.
\begin{Num}
We assume that $D$ is a nonarchimedean and nondyadic local field
(i.e., the characteristic of the residue field of $D$ is $\neq2$).
The Witt group $W\!D$ has $16$ elements,
the group $\widehat S$ of extended square classes $8$, and thus
$I\!I(D,\id)=I^2D$ is cyclic of order $2$ \cite[VI.2.2]{Lam}.
Its nontrivial element is represented by the norm form
of the unique
quaternion division algebra over $D$. Let $S$ denote the group of
square classes of $D$, and 
\[
(-,-)_H:S\times S\rTo\{\pm1\}
\]
the Hilbert symbol \cite[p.~159]{Lam}: $(x,y)_H=-1$ if
$(x,y)_D$ is anisotropic, i.e. the norm form of a quaternion
division algebra. Put $e:I^2D\rTo^\cong\{\pm1\}$, then we 
clearly have
\[
e\circ R\circ stbg(x,y)=(x,y)_H
\] 
for $\SL_2D$.
Thus the reduced Maslov cocycle for $\SL_2D$ is the reduction
of the universal cocycle $stbg$ to $\{\pm1\}$ via the Hilbert symbol.
As in the proof of \ref{CocycleTheorem}, this carries over to $\Sp_{2n}D$.
The following result is partially contained in \cite[p.~104-115]{LV}.
\end{Num}
\begin{Prop}
Let $D$ be a nonarchimedean nondyadic local field. The reduced
Maslov cocycle defines a twofold nontrivial covering 
of $\Sp_{2n}D$ which is determined by the Hilbert symbol
$KSp_2D\rTo\{\pm1\}$. The corresponding covering group
$\widehat{\Sp_{2n}D}$ is a locally compact group; it is the
unique nontrivial twofold covering of $\Sp_{2n}D$ in the category
of locally compact groups.

\proof
Only the topological result remains to be proved. It is shown in
\cite[10.4]{Moore} that in the category of locally compact groups,
$\Sp_{2n}D$ admits a universal central extension
$\widetilde{\Sp_{2n}D}$; the extending
group is the group $\mu(D)$ of all roots of unity in $D$. 
(See \cite{Pra} for a modern account and a much more general result.)
This group $\mu(D)$
is a finite cyclic group \cite[Ch.~II]{Moore} and of even order $2n$,
as it contains the involution
$-1$. The corresponding Steinberg cocycle is given by the
norm residue symbol $KSp_2D\rTo\mu(D)$ \cite[Ch.~II]{Moore}.
But the $n$th power of
the norm residue symbol is the Hilbert symbol.
This shows that $\widehat{\Sp_{2n}D}$ is a continuous
quotient of $\widetilde{\Sp_{2n}D}$.
As the cyclic group $\mu(D)$ has a unique subgroup of index $2$,
the extension is the unique nonsplit twofold topological
extension.
\qed
\end{Prop}
\begin{Num}
\label{UnitarySituation}
Finally, we consider unitary groups over fields. We assume that
$E$ is a field with an automorphism $J\neq\id$ of order $2$;
the fixed field is $D\subseteq E$. We denote the hyperbolic 
unitary group by $\U_{2n}E$; then
$\mathbf{EU}_{2n}E=\SU_{2n}E=\U_{2n}E\cap\SL_{2n}E$ 
\cite[6.4.25,6.4.27]{HOM}.
As we noted in Section \ref{Naturality}, there is
a natural injection $\Phi:\Sp_{2n}D\rInto\SU_{2n}E$ and we have
a commutative diagram
\begin{diagram}[width=5em]
H_2(\Sp_{2n}D)&\rTo^{\Phi_*}&H_2(\SU_{2n}E)\\
\dTo^{[\widetilde\mas_D]}&&\dTo_{[\widetilde\mas_E]}\\
I\!I(D,\id)&\rTo^{W^D_E}&I\!I(E,J).
\end{diagram}
Unfortunately, the Schur multiplier $H_2(\SU_{2n}E)$ seems to be
less understood than its symplectic counterpart. However
it is proved in \cite[2.1,2.5]{Deo} (and in a weaker form in
\cite[6.5.12]{HOM}) that the map
$\Phi_*$ is surjective, so $H_2(\SU_{2n}E)$ is a quotient of
$KSp_2D$. ($\SU_{2n}E$ is the group of $D$-points of a quasisplit
absolutely simple and simply connected algebraic group over $D$, so 
the results from \cite{Deo} apply.) 

The following facts concerning $W^D_E$
were kindly pointed out by W.~Scharlau. Firstly, the map
\[
W^D_E:W^1(D,\id)\rTo W^1(E,J)
\]
is an epimorphism, because every hermitian form can be diagonalized
(even in characteristic $2$ \cite[I.6.2.4]{Knus}) and thus is the
image of a diagonal symmetric bilinear form over $D$.

Assume now that $\mathrm{char}(D)\neq 2$ and $E=D(\sqrt\delta)$.
Passing from a hermitian form $h$ over $E$ to its trace form 
$b_h$ over $D$
\cite[p.~348]{Sch}, we have an monomorphism
$\mathit{trf}:W^1(E,J)\rTo W^1(D,\id)=W\!D$; explicitly,
\[
\mathit{trf}\bra{a_1,\ldots,a_n}=\bra{1,-\delta}\otimes\bra{a_1,\ldots,a_n}.
\]
In particular, $\mathit{trf}\circ W^D_E((x,y)_D)=\bra{1,-\delta}\otimes
\bra{1,-x}\otimes\bra{1,-y}$, and $W^D_E(I^2D)$ is isomorphic
to a subgroup of $I^3D$.
\end{Num}
It follows that the Maslov cocycle for the unitary group over a
nonarchimedean nondyadic local field $E$ vanishes, because $I^3D=0$
\cite[VI.2.15(3)]{Lam}. The case of the complex numbers is more
interesting.
\begin{Num}\textbf{Complex unitary groups}\\
For $E/D=\CC/\RR$ the map
$W^\RR_\CC:W\RR\rTo W^1(\CC,\bar{\phantom{-}})$
and its restriction $I^2\RR\rTo I\!I(\CC,\bar{\phantom{-}})$
is an isomorphism. We use
the standard Lie group notation $\SU_{2n}\CC=\SU(n,n)$ \cite{He} 
(note that multiplication by $i$ transforms skew hermitian into
hermitian matrices).
The maximal compact subgroup is $\mathbf S(\U(n)\times\U(n))$.
As in \ref{RealMaslov}
we compare the classifying space for the discrete group 
(whose homology is the abstract group homology) with the
classifying space $B\SU(n,n)$ for the Lie group.
For $n=1$ we have an isomorphism $\Sp_2\RR=\SU_2\CC$,
whence a big commutative diagram\pagebreak
\begin{diagram}
H_2(\Sp_2\RR)&&\rEq&&H_2(\SU(1,1))&&\rTo^{F_*}&&H_2(B\SU(1,1))\\
\dTo&\rdTo^{[\widetilde\mas_\RR]}&&&\dTo&\rdTo^{[\widetilde\mas_\CC]}
&&&\dTo_\cong&\rdTo^{c_1}_\cong\\
&&I^2\RR &\rLine^{\phi_*}_\cong&\VonH&\rTo&I\!I(\CC,\bar\ )&
\lLine_\cong^{\frac14sig}
&\VonH&\rTo&\ZZ\\
&\ruTo_{[\widetilde\mas_\RR]} &&&&\ruTo_{[\widetilde\mas_\CC]}&&&&
\ruTo_{c_1}^\cong\\
H_2(\Sp_{2n}\RR)&&\rTo&&H_2(\SU(n,n))&&\rTo^{F_*}&&H_2(B\SU(n,n)).
\end{diagram}
\end{Num}
\begin{Prop}
If we identify the first Chern class $c_1$ with the generator
of $H^2(B\SU(n,n))$, it pulls under the forgetful map
$F$ back to the reduced Maslov cocycle for $\SU(n,n)$.
Thus $\widehat{\SU(n,n)}$ is the universal covering group of
the Lie group $\SU(n,n)$.
\qed
\end{Prop}
\begin{Num}
Sharpe \cite{Sharpe} \cite[5.6D*]{HOM} has constructed an exact
sequence
\[
K_2(D)\rTo KU_2^{-1}(D,J)\rTo L_0^1(D,J)\rTo0
\]
The $L$-group $L_0^1(D,J)$ maps onto $I\!I(D,J)$ and we conjecture
that the composite 
\[
KU_2^{-1}(D,J)\rTo I\!I(D,J)
\]
'is' (in most cases) the reduced Maslov cocycle
$\widetilde\mas:H_2(\EU(M))\rTo II(D,J)$. In the symplectic
situation over fields of characteristic $\neq 2$,
this is indeed the case by \ref{CocycleTheorem} and
\cite[5.6.8]{HOM}. However, a proof would certainly require
a different description of the relevant maps than the one in
\cite{Sharpe}.
\end{Num}

\bibliographystyle{plain}

\begin{thebibliography}{WWW}

\bibitem{Ba}
J. Barge,
Cocycle d'Euler et $K\sb 2$, $K$-Theory {\bf 7} (1993),
no.~1, 9--16. MR1220423 (94d:19003)

\bibitem{Bro}
K. S. Brown,
{\it Cohomology of groups},
Corrected reprint of the 1982 original,
Springer, New York, 1994. MR1324339 (96a:20072)

\bibitem{CLM}
S. E. Cappell, R. Lee\ and\ E. Y. Miller,
On the Maslov index,
Comm. Pure Appl. Math. {\bf 47} (1994), no.~2, 121--186.
MR1263126 (95f:57045)

\bibitem{Deo}
Vinay V. Deodhar,
On central extensions of rational points of algebraic groups,
Amer. J. Math. 100 (1978), no. 2, 303--386.
MR489962 (80c:20058)

\bibitem{ES}
S. Eilenberg\ and\ N. Steenrod,
{\it Foundations of algebraic topology}, Princeton Univ. Press, 
Princeton, New Jersey, 1952. MR0050886 (14,398b)

\bibitem{HOM}
A. J. Hahn\ and\ O. T. O'Meara,
{\it The classical groups and $K$-theory},
Springer, Berlin, 1989. MR1007302 (90i:20002)

\bibitem{He}
S. Helgason,
{\it Differential geometry, Lie groups, and symmetric spaces},
Corrected reprint of the 1978 original,
Amer. Math. Soc., Providence, RI, 2001. MR1834454 (2002b:53081)

\bibitem{Knarr}
N. Knarr,
Projectivities of generalized polygons,
Ars Combin. {\bf 25} (1988), B, 265--275.
MR0942482 (89e:20008)

\bibitem{Knus}
M.-A. Knus,
{\it Quadratic and Hermitian forms over rings},
Springer, Berlin, 1991. MR1096299 (92i:11039) 

\bibitem{LK}
L. Kramer,
Buildings and classical groups,
in {\it Tits buildings and the model theory of groups
(W\"urzburg, 2000)}, 59--101, Cambridge Univ. Press,
Cambridge. MR2018382 (2005b:20058)

\bibitem{Kub}
T. Kubota,
Topological covering of ${\rm SL}(2)$ over a local field,
J. Math. Soc. Japan {\bf 19} (1967), 114--121.
MR0204422 (34 \#4264)

\bibitem{Lam}
T. Y. Lam,
{\it Introduction to quadratic forms over fields},
Amer. Math. Soc., Providence, RI, 2005. 
MR2104929 (2005h:11075)

\bibitem{LV}
G. Lion\ and\ M. Vergne,
{\it The Weil representation, Maslov index and theta series},
Progr. Math., 6, Birkh\"auser, Boston, Mass., 1980.
MR0573448 (81j:58075)

\bibitem{McL}
S. Mac Lane,
{\it Homology}, Reprint of the 1975 edition,
Springer, Berlin, 1995. MR1344215 (96d:18001)

\bibitem{Mat}
H. Matsumoto,
Sur les sous-groupes arithm\'etiques des groupes semi-simples
d\'eploy\'es,
Ann. Sci. \'Ecole Norm. Sup. (4) {\bf 2} (1969), 1--62.
MR0240214 (39 \#1566)

\bibitem{Maz}
A. Mazzoleni,
Partially defined cocycles and the Maslov index for a local ring,
Ann. Inst. Fourier (Grenoble) {\bf 54} (2004),
no.~4, 875--885. MR2111015 (2005h:20119)

\bibitem{Moore}
C. C. Moore,
Group extensions of $p$-adic and adelic linear groups,
Inst. Hautes \'Etudes Sci. Publ. Math. No. 35 (1968), 157--222. 
MR0244258 (39 \#5575)

\bibitem{NO}
K.-H. Neeb\ and\ B. \O rsted,
A topological Maslov index for 3-graded Lie groups,
J. Funct. Anal. {\bf 233} (2006), no.~2, 426--477.
MR2214583 (2007h:53127)

\bibitem{Ne}
Ya. Nekovar,
Maslov index and Clifford algebras,
Funktsional. Anal. i Prilozhen. {\bf 24} (1990), no. 3, 36--44, 96;
translation in Funct. Anal. Appl. {\bf 24} (1990), no.~3, 196--204 (1991).
MR1082029 (92b:11024)

\bibitem{No}
M. V. Nori,
The universal property of the Maslov index,
J. Ramanujan Math. Soc. {\bf 13} (1998), no.~2, 111--124.
MR1666437 (2000e:11039

\bibitem{PPS}
R. Parimala, R. Preeti\ and\ R. Sridharan,
Maslov index and a central extension of the symplectic group,
$K$-Theory {\bf 19} (2000), no.~1, 29--45. MR1740881 (2001c:11053a)

R. Parimala, R. Preeti\ and\ R. Sridharan,
Errata: ``Maslov index and a central extension of the symplectic group'',
$K$-Theory {\bf 19} (2000), no.~4, 403. MR1763935 (2001c:11053b)

\bibitem{Pra}
G. Prasad,
Deligne's topological central extension is universal,
Adv. Math. {\bf 181} (2004), no.~1, 160--164.
MR2020658 (2004k:20097)

\bibitem{Sch}
W. Scharlau,
{\it Quadratic and Hermitian forms},
Springer, Berlin, 1985. MR0770063 (86k:11022)

\bibitem{Schwarze}
R. Schwarze,
{\it Group extensions and the Maslov index},
Diploma Thesis, Bielefeld, 2008.

\bibitem{Sharpe}
R. W. Sharpe,
On the structure of the unitary Steinberg group,
Ann. of Math. (2) {\bf 96} (1972), 444--479.
MR0320076 (47 \#8617)

\bibitem{Spa}
E. H. Spanier,
{\it Algebraic topology},
Corrected reprint, Springer, New York, 1981. MR0666554 (83i:55001)

\bibitem{Ti} J. Tits,
{\it Buildings of spherical type and finite BN-pairs},
Lecture Notes in Math., 386, Springer, Berlin, 1974. 
MR0470099 (57 \#9866)

\bibitem{Tu}
V. G. Turaev,
A cocycle of the symplectic first Chern class and Maslov indices,
Funktsional. Anal. i Prilozhen. {\bf 18} (1984), no.~1, 43--48.
MR0739088 (85m:58191)

\bibitem{Weil}
A. Weil,
{\it Basic number theory},
Reprint of the second (1973) edition, Springer, Berlin, 1995.
MR1344916 (96c:11002)

\end{thebibliography}

\bigskip
\raggedright
Linus Kramer\\
Katrin Tent\\
Mathematisches Institut, 
Universit\"at M\"unster,
Einsteinstr. 62,
48149 M\"unster,
Germany\\
\makeatletter
e-mail: {\tt linus.kramer{@}uni-muenster.de}\\
e-mail: {\tt tent{@}uni-muenster.de}

\end{document}